\documentclass[12pt]{article}

\usepackage{geometry}
\usepackage{amsmath}
\usepackage{amssymb}
\usepackage{bbm}

\usepackage{graphicx}
\graphicspath{ {./figures/} }

\usepackage{enumitem}

\usepackage{xcolor}

\usepackage{algorithm}
\usepackage{caption}
\usepackage[algo2e,linesnumbered,ruled]{algorithm2e}

\SetCommentSty{mycommfont}
\let\oldnl\nl
\newcommand{\nonl}{\renewcommand{\nl}{\let\nl\oldnl}}

\newtheorem{definition}{Definition}[section]
\newtheorem{lemma}{Lemma}[section]
\newtheorem{proposition}{Proposition}[section]
\newtheorem{theorem}{Theorem}[section]
\newtheorem{remark}{Remark}[section]
\newtheorem{example}{Example}[section]
\newtheorem{assumption}{Assumption}[section]

\usepackage[pagebackref]{hyperref}       
\hypersetup{
	colorlinks=true,
	linkcolor=blue,
	filecolor=blue,
	citecolor=blue,      
	urlcolor=blue,
}
\renewcommand*{\backref}[1]{}
\renewcommand*{\backrefalt}[4]{%
	\ifcase #1 Not cited.%
	\or        Cited on page~#2.%
	\else      Cited on pages~#2.%
	\fi}

\usepackage{cleveref}
\crefname{figure}{Figure}{Figures}
\crefname{section}{Section}{Sections}
\crefname{lemma}{Lemma}{Lemmas}
\crefname{proposition}{Proposition}{Propositions}
\crefname{theorem}{Theorem}{Theorem}
\crefname{definition}{Definition}{Definitions}
\crefname{remark}{Remark}{Remarks}
\crefname{algorithm}{Algorithm}{Algorithms}
\crefname{appendix}{Appendix}{Appendices}
\crefname{example}{Example}{Examples}
\crefname{assumption}{Assumption}{Assumptions}

\usepackage{crossreftools}
\makeatletter
\renewcommand\paragraph{%
	\@startsection{paragraph}
	{4}
	{\z@}
	{3.25ex \@plus1ex \@minus.2ex}
	{-1em}
	{\normalfont\normalsize\bfseries\maybe@addperiod}%
}
\newcommand{\maybe@addperiod}[1]{%
	#1\@addpunct{.}%
}
\makeatother

\DeclareMathOperator*{\argmax}{arg\,max}

\newcommand{\succsucc}{\succ\mathrel{\mkern-5mu}\succ}

\usepackage[parfill]{parskip}

\usepackage[nottoc,numbib]{tocbibind}

\usepackage[font=footnotesize,labelfont=bf, textfont=it]{caption}

\makeatletter
\def\@fnsymbol#1{\ensuremath{\ifcase#1\or \dagger\or \ddagger\or
		\mathsection\or \mathparagraph\or \|\or **\or \dagger\dagger
		\or \ddagger\ddagger \else\@ctrerr\fi}}
\makeatother

\usepackage[parfill]{parskip}

\title{Multi-objective optimisation via the R2 utilities}

\author{
	Ben Tu\thanks{Imperial College London, United Kingdom}
	\and Nikolas Kantas\footnotemark[1]
	\and Robert M. Lee\thanks{BASF SE, Germany}
	\and Behrang Shafei\footnotemark[2]
}
\date{}

\begin{document}
	
\maketitle

\begin{abstract}
	The goal of multi-objective optimisation is to identify a collection of points which describe the best possible trade-offs among the multiple objectives. In order to solve this vector-valued optimisation problem, practitioners often appeal to the use of scalarisation functions in order to transform the multi-objective problem into a collection of single-objective problems. This set of scalarised problems can then be solved using traditional single-objective optimisation techniques. In this paper, we formalise this convention into a general mathematical framework. We show how this strategy effectively recasts the original multi-objective optimisation problem into a single-objective optimisation problem defined over sets. An appropriate class of objective functions for this new problem is that of the R2 utilities, which are utility functions that are defined as a weighted integral over the scalarised optimisation problems. As part of our work, we show that these utilities are monotone and submodular set functions which can be optimised effectively using greedy optimisation algorithms. We then analyse the performance of these greedy algorithms both theoretically and empirically. Our analysis largely focusses on Bayesian optimisation, which is a popular probabilistic framework for black-box optimisation.
\end{abstract}

\section{Introduction}
Decisions that are made in the real-world lead to a variety of consequences that have to be considered beforehand. It is rare that a decision leads only to positive changes that would benefit everyone involved. In many cases, there is a trade-off that has to be struck which tries to balance between many different criteria. Multi-objective optimisation formalises this decision making problem mathematically as a vector-valued optimisation problem. The solution to this problem is a set of decisions, whose objective values lead to the best possible trade-offs among the multiple criteria. Equipped with this solution set, a decision maker can then make a more informed decision which takes into account all the possible trade-offs.

Many efforts over the past few decades have been spent developing the theory and algorithms that are used to solve this multi-objective optimisation problem \cite{miettinen1998,ehrgott2005}. The aim of this paper is to consolidate the core ideas of these approaches into a mathematical framework that will support the future development of algorithms and the theory behind them. Notably, we study the family of R2 utility functions, which connects three different fields of optimisation: multi-objective optimisation, single-objective optimisation and set optimisation. We support these connections with general theoretical results and illustrative examples that provide intuition. Whilst these results are not surprising and are already known for a specific instance of R2 utility, namely, the hypervolume indicator \cite{zitzler1998ppsn,ulrich2012lio}, to the best of our knowledge they have not been presented in the most general form and are absent from the literature. Furthermore, these clear connections can also be exploited in order to design more principled algorithms. To illustrate this, we present a motivating example of how these ideas can be applied for Bayesian optimisation, which is a popular black-box optimisation strategy that has gained a lot of interest and success over the recent decade. We also provide general performance bounds for this case and detailed numerical examples.

\subsection{Structure of the paper} The remainder of the paper is organised as follows: In \cref{sec:preliminaries}, we define the multi-objective optimisation problem and introduce the two main philosophies that are often used to solve this problem: the scalarisation perspective and the utility perspective. In \cref{sec:r2_utility}, we unify these two methodologies by introducing the R2 utilities, which form a family of utility functions that are defined using scalarisation functions. We then show that the family of R2 utilities satisfies many desirable properties, which makes them a sensible criterion to optimise. In \cref{sec:optimisation_problem}, we give a discussion on how we can solve the multi-objective optimisation problem using R2 utilities. We focus our attention on greedy optimisation strategies because they come with theoretical guarantees. In particular, we prove a general performance bound which is satisfied by any greedy Bayesian optimisation algorithm based on R2 utilities. In \cref{sec:experiments}, we evaluate these greedy Bayesian optimisation strategies empirically in light of these new theoretical results. Finally, in \cref{sec:future_work}, we conclude the paper with a summary and discussion of future work. In \cref{app:proofs}, we give the proofs of the main results. 
\section{Preliminaries}
\label{sec:preliminaries}
Consider a vector-valued function $f: \mathbb{X} \rightarrow \mathbb{R}^M$ defined over a $D$-dimensional space of feasible inputs $\mathbb{X} \subseteq \mathbb{R}^D$. The multi-objective maximisation problem is denoted by the equation
\begin{equation}
	\max_{\mathbf{x} \in \mathbb{X}} f(\mathbf{x}),
	\label{eqn:multi-objective_problem}
\end{equation}
where the maximum is defined via the Pareto partial ordering relation.
\begin{definition}
	[Pareto domination] The weak, strict, and strong Pareto domination are denoted by the binary relations $\succeq, \succ$ and $\succsucc$, respectively. We say that a vector $\mathbf{y} \in \mathbb{R}^M$ weakly, strictly, or strongly Pareto dominates another vector $\mathbf{y'} \in \mathbb{R}^M$, if
	\begin{align*}
		\mathbf{y} \succeq \mathbf{y'} &\iff \mathbf{y} - \mathbf{y'} \in \mathbb{R}_{\geq 0}^M,
		\\
		\mathbf{y} \succ \mathbf{y'} &\iff \mathbf{y} - \mathbf{y'} \in \mathbb{R}_{\geq 0}^M \setminus \{\mathbf{0}_M\},
		\\
		\mathbf{y} \succsucc \mathbf{y'} &\iff \mathbf{y} - \mathbf{y'} \in \mathbb{R}_{> 0}^M,
	\end{align*}
	respectively, where $\mathbf{0}_M \in \mathbb{R}^M$ denotes the $M$-dimensional vector of zeros.
\end{definition}
\begin{definition}
	[Pareto optimality] Consider a function $f: \mathbb{X} \rightarrow \mathbb{R}^M$, we say an input $\mathbf{x} \in \mathbb{X}$ is weakly or strictly Pareto optimal if the objective vector $f(\mathbf{x})$ is not strongly or strictly dominated, respectively, by any other objective vector $f(\mathbf{x}')$ with $\mathbf{x}' \in \mathbb{X} \setminus \{\mathbf{x}\}$. 
\end{definition}

The convention in multi-objective optimisation is to target the strictly Pareto optimal points. The collection of strictly Pareto optimal inputs is called the Pareto set, $\mathbb{X}^* = \argmax_{\mathbf{x} \in \mathbb{X}} f(\mathbf{x}) \subseteq \mathbb{X}$, whilst the corresponding image is called the Pareto front, $\mathbb{Y}^* = f(\mathbb{X}^*) = \max_{\mathbf{x} \in \mathbb{X}} f(\mathbf{x}) \subset \mathbb{R}^M$. 

\paragraph{Set Pareto domination} Practitioners are often only interested in identifying a discrete approximation of the Pareto front. For this reason, we will consider working in the space of finite sets: $\mathbb{B}(\mathbb{R}^M) = \{Y\subseteq \mathbb{R}^M: |Y| < \infty\}$. In practice, we are interested in making comparisons between different Pareto front approximations. To that end, we will now extend the Pareto partial ordering to be defined over finite sets. Informally, we will say a set $A \in \mathbb{B}(\mathbb{R}^M)$ dominates another set $B \in \mathbb{B}(\mathbb{R}^M)$ if the region dominated by $A$ contains the region dominated by $B$.

\begin{definition}
	[Dominated region] For a set of vectors $Y \in \mathbb{B}(\mathbb{R}^M)$, the weak dominated region is defined as the collection of vectors which is weakly dominated by at least one vector in this set, that is $\mathbb{D}_{\preceq}(Y) = \bigcup_{\mathbf{y} \in Y} \{\mathbf{a} \in \mathbb{R}^M: \mathbf{a} \preceq \mathbf{y}\}$.
\end{definition}

\begin{definition}
	[Set Pareto domination] We say a set of vectors $A \in \mathbb{B}(\mathbb{R}^M)$ weakly or strictly Pareto dominates another set of vectors $B \in \mathbb{B}(\mathbb{R}^M)$, if
	\begin{align*}
		A \succeq B &\iff \mathbb{D}_{\preceq}(A) \supseteq \mathbb{D}_{\preceq}(B),
		\\
		A \succ B &\iff \mathbb{D}_{\preceq}(A) \supset \mathbb{D}_{\preceq}(B),
	\end{align*}
	respectively. 
	\label{def:set_pareto_domination}
\end{definition}

\begin{remark}
	The above definition of set domination might differ slightly from alternative definitions given in the literature \cite{hansen1998trimm}. Nevertheless, the core idea remains the same: a set dominates another if for any element in the latter set, we can always find an element in the former set which dominates it. 
\end{remark}

\paragraph{Reformulating the multi-objective optimisation problem} Designing algorithms to solve the multi-objective optimisation \eqref{eqn:multi-objective_problem} directly is challenging because the Pareto partial ordering means that not all vectors or sets of vectors are comparable with each other. To address this limitation, decision maker's often rely on the use of scalarisation functions, $s: \mathbb{R}^M \rightarrow \mathbb{R}$, or utility functions, $U: \mathbb{B}(\mathbb{R}^M) \rightarrow \mathbb{R}$, in order make comparisons between vectors or sets of vectors, respectively. In what follows, we will review these two classes of functions and show how they effectively reformulate the multi-objective optimisation problem as a different optimisation problem. We will then unify these two perspectives in \cref{sec:r2_utility}, by showing how it is possible to construct a sensible utility function from a collection of scalarisation functions.
\subsection{Scalarisation perspective}
\label{sec:scalarisation_perspective}
A well-known approach for solving the multi-objective optimisation problem \eqref{eqn:multi-objective_problem} is to transform it into a collection of single-objective optimisation problems. These problems can then be jointly solved using classical techniques from single-objective optimisation. We refer to this reformulation as the scalarisation perspective of multi-objective optimisation because we typically rely on the use of scalarisation functions in order to establish these single-objective problems. This perspective is also sometimes referred to as the decomposition-based approach to multi-objective optimisation \cite{zhou2011saec,li2015acs,trivedi2017itec}.

\paragraph{Scalarised optimisation problems} Formally, in the scalarisation perspective, we are interested in solving a collection of single-objective optimisation problems obtained applying a family of scalarisation functions: $\{s_{\boldsymbol{\theta}}: \mathbb{R}^M \rightarrow \mathbb{R}: \boldsymbol{\theta} \in \Theta\}$. In particular, we want to solve the scalarised optimisation problem,
\begin{equation}
	\max_{\mathbf{x} \in \mathbb{X}} s_{\boldsymbol{\theta}}(f(\mathbf{x}))
	\label{eqn:scalarised_problem}
\end{equation}
jointly for all scalarisation parameters $\boldsymbol{\theta} \in \Theta$. Intuitively, each scalarised problem corresponds to one point that we are interested in targetting. If the number of scalarisation parameters is finite, $|\Theta| < \infty$, then we are only targetting a finite number of points. In contrast, if we had a separate scalarised problem for each strictly Pareto optimal point, then jointly solving this set of problems would be equivalent to solving the original multi-objective optimisation problem.

Ideally, we want the solutions to the scalarised optimisation problems to be Pareto optimal in some sense. This property turns out to be satisfied when the scalarisation function preserves the Pareto partial ordering. Specifically, whenever one vector dominates another, we want the scalarised value to also be greater. This notion of preserving the Pareto partial ordering can be interpreted as a notion of monotonicity.
\begin{definition}
	[Monotonicity] A scalarisation function $s: \mathbb{R}^M \rightarrow \mathbb{R}$ is monotonically increasing, strictly monotonically increasing, or strongly monotonically increasing if
	\begin{align*}
		\mathbf{y} \succeq \mathbf{y}' &\implies s(\mathbf{y}) \geq s(\mathbf{y}'),
		\\
		\mathbf{y} \succ \mathbf{y}' &\implies s(\mathbf{y}) > s(\mathbf{y}'),
		\\
		\mathbf{y} \succsucc \mathbf{y}' &\implies s(\mathbf{y}) > s(\mathbf{y}'),
	\end{align*}
	for any vectors $\mathbf{y}, \mathbf{y}' \in \mathbb{R}^M$, respectively. Similarly, the function is monotonically decreasing, strictly monotonically decreasing, or strongly monotonically decreasing if the reverse inequalities holds.
\end{definition}
A well-known result in multi-objective optimisation states that whenever the scalarisation function is strictly or strongly monotonically increasing, then the solution set is Pareto optimal \cite[Part 2, Theorem 3.5.4]{miettinen1998}---this result is summarised in \cref{prop:monotonic_implies_optimal} and proved in \cref{app:proofs:prop:monotonic_implies_optimal}.
\begin{proposition}
	[Monotonicity implies optimality] Consider an objective function $f:\mathbb{X} \rightarrow \mathbb{R}^M$ and a scalarisation function $s: \mathbb{R}^M \rightarrow \mathbb{R}$. If the scalarisation function is strictly or strongly monotonically increasing over the feasible objective space, then the solutions in $X^*_s = \argmax_{\mathbf{x} \in \mathbb{X}} s(f(\mathbf{x})) \subseteq \mathbb{X}$ are strictly or weakly Pareto optimal, respectively.
	\label[proposition]{prop:monotonic_implies_optimal}
\end{proposition}
\begin{figure}
	\includegraphics[width=1\linewidth]{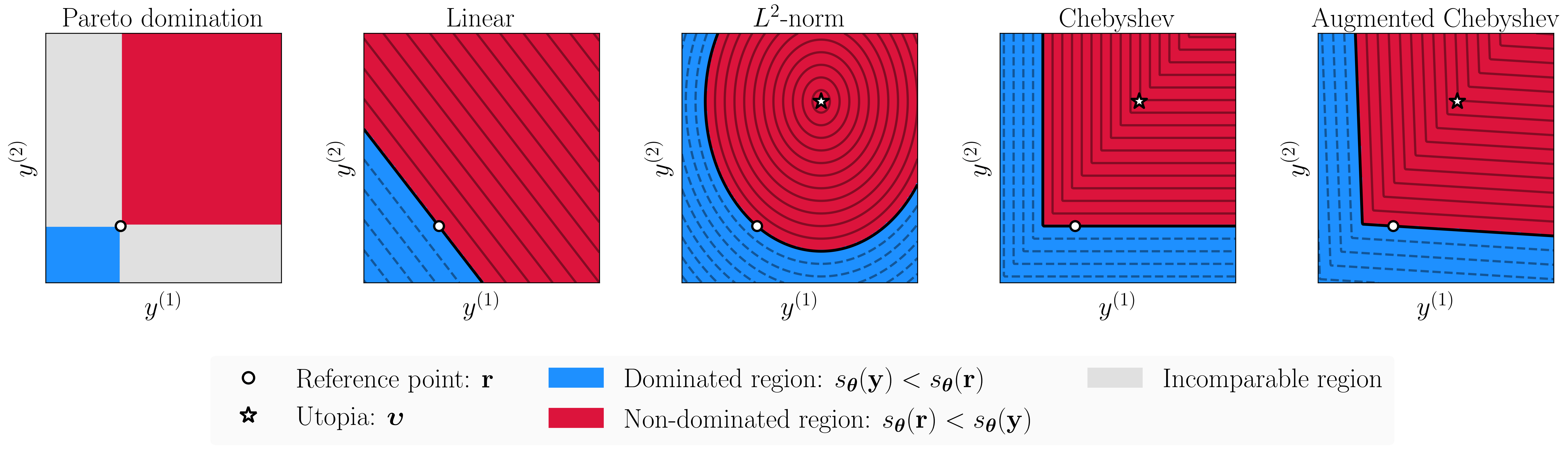}
	\caption{A comparison of the domination regions based on the standard Pareto partial ordering over vectors and some popular scalarisation functions in the two-objective setting.}
	\label{fig:scalarisation_domination}
\end{figure}
\begin{example}
	[Popular scalarisation functions] There are many possible scalarisation functions that we can use in practice \cite{marler2004smo}. Below, we give a few examples of some popular scalarisation functions. In particular, we consider the linear, $L^p$-norm (Lp), Chebyshev (Chb) and augmented Chebyshev (AugChb) scalarisation functions:
	\begin{align}
		s^{\text{Linear}}_{\mathbf{w}}(\mathbf{y}) 
		&= \sum_{m=1}^M w^{(m)} y^{(m)},
		\label{eqn:linear_scalarisation}
		\\
		s^{\text{Lp}}_{(\boldsymbol{\upsilon}, \mathbf{w})}(\mathbf{y}) 
		&= - \biggl(\sum_{m=1}^M |w^{(m)} (\upsilon^{(m)} - y^{(m)})|^p \biggr)^{1/p},
		\label{eqn:lp_scalarisation}
		\\
		s^{\text{Chb}}_{(\boldsymbol{\upsilon}, \mathbf{w})}(\mathbf{y}) 
		&= - \max_{m = 1, \dots, M} w^{(m)}(\upsilon^{(m)} - y^{(m)}),
		\label{eqn:chebyshev_scalarisation}
		\\
		s^{\text{AugChb}}_{(\boldsymbol{\upsilon}, \mathbf{w}, \gamma)}(\mathbf{y}) 
		&= s^{\text{Chb}}_{(\boldsymbol{\upsilon}, \mathbf{w})}(\mathbf{y}) 
		+  \gamma s^{\text{Linear}}_{\mathbf{w}}(\mathbf{y}),
		\label{eqn:augmented_chebyshev_scalarisation}
	\end{align}
	respectively, where $\mathbf{y} \in \mathbb{R}^M$ is an objective vector, $\mathbf{w} \in \Delta^{M-1}:= \{\mathbf{y} \in \mathbb{R}^M_{\geq 0}: ||\mathbf{y}||_{L^1} = 1\}$ is a weight vector lying in the non-negative $M$-dimensional simplex, $\boldsymbol{\upsilon} \in \mathbb{R}^M$ is an ideal\footnote{An ideal or utopia point is a vector whose objective values are desirable and conversely, the disagreement or nadir point is a vector comprised of undesirable objective values. Both of these vectors are typically set according to the decision maker's preferences \cite{wagner2013emo}.} reference point, $p \geq 1$ is a constant controlling the $L^p$-norm, and $\gamma \geq 0$ is a penalty parameter.  In \cref{fig:scalarisation_domination}, we illustrate the contours of these scalarisation functions in the bi-objective setting, for one particular choice of scalarisation parameter. 
	\par
	Notably, the linear, Chebyshev, and augmented Chebyshev scalarisation functions are strongly monotonically increasing over the whole objective space and therefore they satisfy the result of \cref{prop:monotonic_implies_optimal}. On the other hand, the Lp scalarisation function is only strongly monotonically increasing in the space dominated by the reference point and therefore it only satisfies this result when the reference point is set accordingly. 
	\label{eg:scalarisation_functions}
\end{example}

\begin{remark}
	[Chebyshev scalarisation] \cref{prop:monotonic_implies_optimal} ensures that some Pareto optimal solutions can be targetted when we vary the scalarisation parameter. It does not give any guarantees on targetting all of the Pareto optimal solutions. These type of results are typically only possible for explicit choices of scalarisation functions. For example, it is known that we can target all the weakly Pareto optimal solutions, which are strictly Pareto comparable with the reference vector $\boldsymbol{\upsilon} \in \mathbb{R}^M$, by varying the weight parameter $\mathbf{w} \in \Delta^{M-1}$ in the Chebyshev scalarisation \eqref{eqn:chebyshev_scalarisation} function \cite[Part 2, Theorem 3.4.5]{miettinen1998}. 
	\label{rem:chebyshev}
\end{remark}

\begin{remark}
	[Objective transformations] Scalarisation functions are naturally sensitive to the scales of the objectives. Thus, it is common for practitioners to apply a transformation function, $\tau: \mathbb{R}^M \rightarrow \mathbb{R}^M$, to the objective vector before applying the scalarisation function $s(\tau(\mathbf{y})) \in \mathbb{R}$. For instance, one might consider applying a linear transformation, which normalises the feasible objective values to the unit hypercube, $[0, 1]^M$, before computing any of the distance-based scalarisation functions in \cref{eg:scalarisation_functions}. In this paper, we will implicitly assume that any such transformation has already been included in the definition of the scalarisation function.
\end{remark}

\subsection{Utility perspective}
\label{sec:utility_perspective}
One of the primary challenges in multi-objective optimisation concerns the quantification of the quality of an approximate Pareto front. As observed by Zitzler et al. \cite{zitzler2008moiaea}, the set Pareto domination is only a partial ordering relation, which means that we are not always able to determine whether one approximation set is better than another. Therefore, the burden often falls onto the decision maker to characterise whether one set is more preferable than another, even when the sets are not Pareto comparable. To that end, decision maker's often rely on utility functions, $U: \mathbb{B}(\mathbb{R}^M) \rightarrow \mathbb{R}$, in order to determine the quality of an approximate Pareto front. This naturally leads to the utility perspective of multi-objective optimisation, which considers the problem of directly optimising this utility function. This perspective is also sometimes referred to as the indicator-based approach to multi-objective optimisation \cite{zhou2011saec,li2015acs,trivedi2017itec}.

\paragraph{Utility optimisation problem} Formally, in the utility perspective, we are interested in optimising the utility function
\begin{equation}
	\max_{X \subseteq \mathbb{X}, |X| \leq P} U(f(X))
	\label{eqn:utility_problem}
\end{equation}
for some number $P > 0$. We have included a cardinality constraint into this set optimisation problem in order to emphasise the fact that we are typically interested in building a discrete approximation of the Pareto front. Naturally, we could lift this constraint if an abundance of computational power was available.

Ideally, the utility function should be designed to preserve the Pareto partial ordering over sets. Similar to the notion of monotonicity for scalarisation functions, we introduce the notion of compliancy for utility functions. 
\begin{definition}
	[Pareto compliancy] A utility function $U: \mathbb{B}(\mathbb{R}^M) \rightarrow \mathbb{R}$ is weakly or strictly Pareto compliant if
	\begin{align*}
		A \succeq B &\implies U(A) \geq U(B), \\
		A \succ B &\implies U(A) > U(B),
	\end{align*}
	for any sets of vectors $A, B \in \mathbb{B}(\mathbb{R}^M)$, respectively.
	\label{def:pareto_compliant}
\end{definition}

The Pareto compliancy property implies that the utility optimisation problem \eqref{eqn:utility_problem} can be solved by using some subset of the Pareto set. When the Pareto front is finite\footnote{Extending this result to the infinite setting is non-trivial. Specifically, we would have to deal with some measurability issues involved with extending \cref{def:set_pareto_domination}.}, the strict Pareto compliancy property leads to a stronger result which states that the smallest set which maximises the utility is the Pareto front; that is, $\max_{X \in \mathbb{X}, |X| \leq P} U(f(X)) = U(f(\mathbb{X}^*)) = U(\mathbb{Y}^*)$ for any $P\geq |\mathbb{Y}^*|$. 
\begin{figure}
	\includegraphics[width=0.24\linewidth]{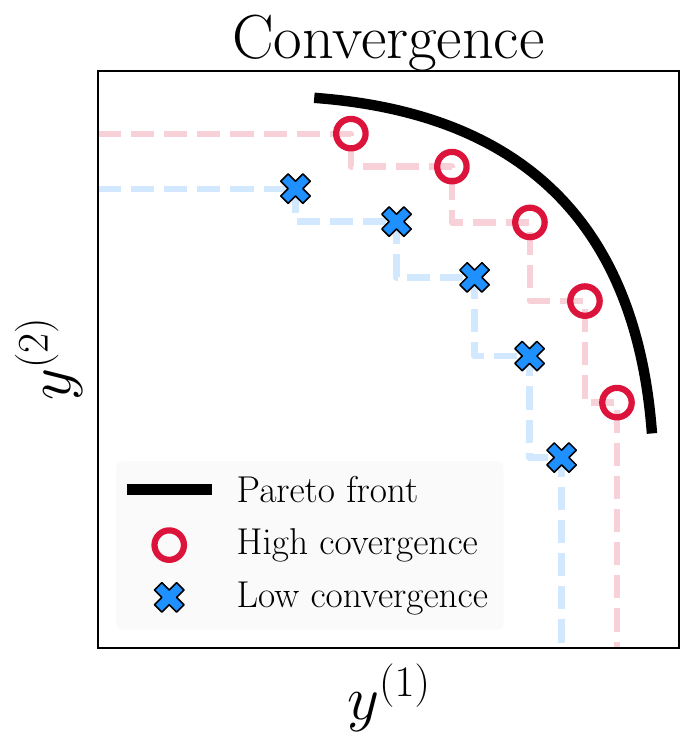}
	\includegraphics[width=0.24\linewidth]{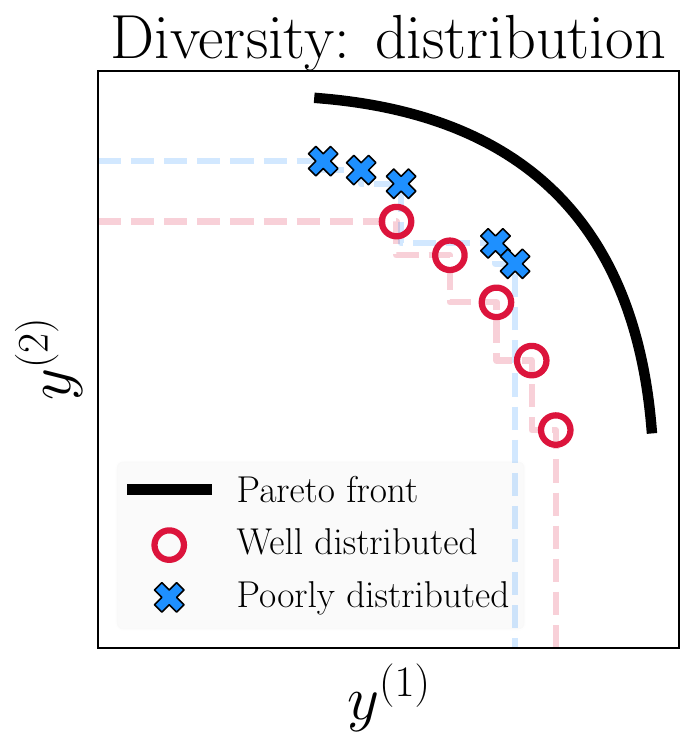}
	\includegraphics[width=0.24\linewidth]{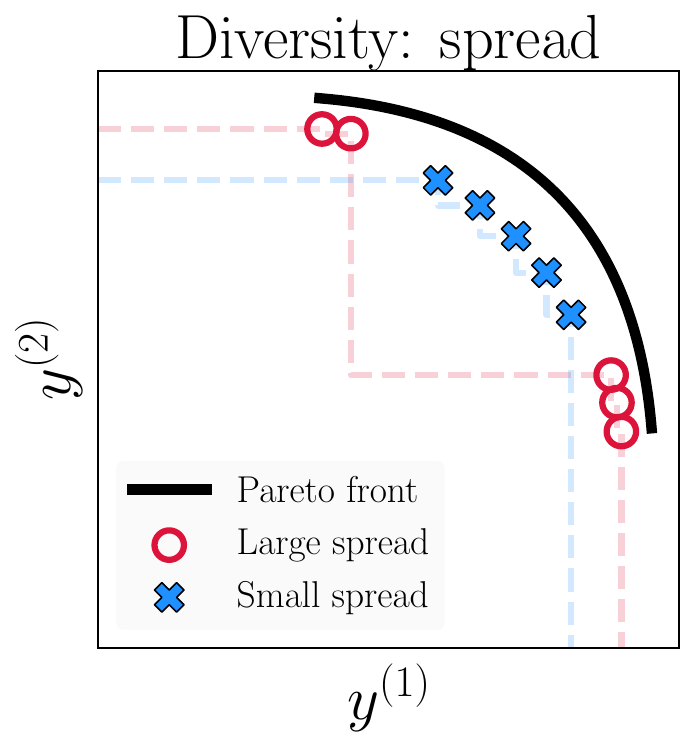}
	\includegraphics[width=0.24\linewidth]{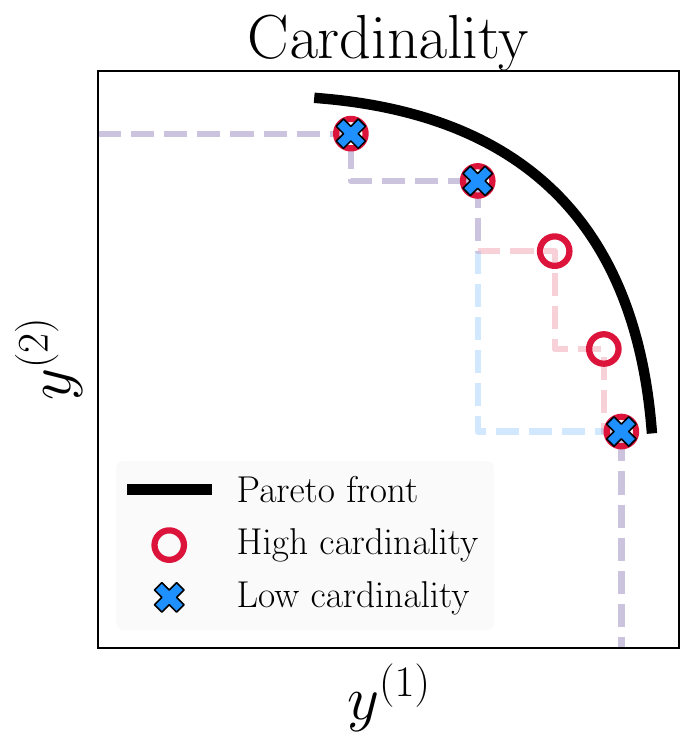}
	\caption{A comparison of the aspects that determine the quality of a Pareto front approximation.}
	\label{fig:three_quals}
\end{figure}

\paragraph{Qualitative aspects}
In the multi-objective literature, the utility function is often referred to as a performance metric or a performance indicator. Following Riquelme, von Lücken, and Baran \cite{riquelme20152laccc}, there are three core qualitative aspects that a multi-objective performance metric inherently tries to quantify.

\begin{enumerate}
	\item \textbf{Convergence.} This aspect is concerned with the accuracy of the approximation, namely, how ``close'' the approximate front is to the true Pareto front. This aspect requires the specification of some notion of distance in the objective space, which can be especially challenging  when the objectives describe quantities that are defined on different scales. 
	\item \textbf{Diversity.} This aspect is concerned with the distribution and spread of the points in the approximate front. The distribution is related to the relative distance among the points, whereas the spread considers the range of the objectives that are covered by the points. As illustrated in \cref{fig:three_quals}, a set can have a good distribution, but a poor spread and vice versa. 
	\item \textbf{Cardinality.} This aspect refers to the number of points in the approximate front. Naturally, a larger number of points is typically preferred.  
\end{enumerate}

In \cref{fig:three_quals}, we illustrate these three qualities in the two-objective setting. As observed in this example, an approximation can do well on one of these aspects and poorly on another. In essence, all performance metrics in multi-objective optimisation work by quantifying some or all of these aspects and then striking a suitable balance between them. On the surface it appears that there are many different ways to design a performance metric which strikes a balance among these three qualitative properties. For instance, a recent survey \cite{audet2021ejoor} identified over fifty different performance metrics that have been documented in the literature. In this paper, we will focus our attention on the \emph{R2 utilities}, which form a general class of utility functions that are inspired by the R2 metric \cite{hansen1998trimm}. This family turns out to possess many desirable properties and contains many of the most often used performance metrics as special cases---more details are given in \cref{sec:r2_utility}. 

\paragraph{Other useful properties} We will now define two other desirable properties which we will need later on. The first definition is the monotone property, which states that a utility function is non-decreasing in the sense that adding more points to a set does not decrease its utility. The second definition is the diminishing returns property, which states that the gain in utility is less for larger sets. In this paper, we will say a set function, $U: \mathbb{B}(\Omega) \rightarrow \mathbb{R}$, is \emph{submodular} if it satisfies the diminishing returns property. We will consider the setting where the ground set is the space of vectors: $\Omega = \mathbb{R}^M$. Note that this definition of submodularity differs subtly from the conventional setting \cite{bach2013fiml, krause2014tpathp}, which additionally assumes that the ground set is finite.
\begin{definition}
	[Monotone] For a ground set $\Omega$, the set function $U: \mathbb{B}(\Omega) \rightarrow \mathbb{R}$ is monotone if for any sets $A, B \in \mathbb{B}(\Omega)$ with $A \subseteq B$ we have that $U(A) \leq U(B)$.
	\label[definition]{def:monotone}
\end{definition}
\begin{definition}
	[Diminishing returns] For a ground set $\Omega$, the set function $U: \mathbb{B}(\Omega) \rightarrow \mathbb{R}$ satisfies the diminishing property if and only if for any sets $A, B \in \mathbb{B}(\Omega)$ with $A \subseteq B$ and $\mathbf{c} \in \Omega$ we have that $U(A \cup \{\mathbf{c}\}) - U(A) \geq U(B \cup \{\mathbf{c}\}) - U(B)$.
	\label[definition]{def:diminishing_returns}
\end{definition}

\section{R2 Utilities} 
\label{sec:r2_utility}
The family of R2 utilities offers a natural way to connect between the scalarisation and utility perspectives of multi-objective optimisation, by defining a utility function based on a family of scalarisation functions. In particular, given a family of scalarisation functions $\{s_{\boldsymbol{\theta}}: \mathbb{R}^M \rightarrow \mathbb{R}: \boldsymbol{\theta} \in \Theta\}$ and a probability density $p(\boldsymbol{\theta}) \geq 0$, the corresponding R2 utility function is defined as the average maximum scalarised value
\begin{equation}
	U(Y) 
	= \mathbb{E}_{p(\boldsymbol{\theta})}[S_{\boldsymbol{\theta}}(Y)]
	= \mathbb{E}_{p(\boldsymbol{\theta})}\biggl[
	\max_{\mathbf{y} \in Y} s_{\boldsymbol{\theta}}(\mathbf{y})
	\biggr]
	\label{eqn:r2_utility}
\end{equation}
for some set of vectors $Y \in \mathbb{B}(\mathbb{R}^M)$, where $\mathbb{E}_{p(\boldsymbol{\theta})}[\cdot]$ denotes the expectation with respect to the density $p(\boldsymbol{\theta})$. Intuitively, the maximum scalarised value, described by the set function, $S_{\boldsymbol{\theta}}: \mathbb{B}(\mathbb{R}^M) \rightarrow \mathbb{R}$,
\begin{equation*}
	S_{\boldsymbol{\theta}}(Y) = \max_{\mathbf{y} \in Y} s_{\boldsymbol{\theta}}(\mathbf{y}),
\end{equation*}
evaluates the quality of targetting one point corresponding to the scalarised optimisation problem \eqref{eqn:scalarised_problem}. By taking a weighted average over these problems, the R2 utility quantifies some notion of average quality when we are interested in targetting multiple points. To keep the discussion focussed, we will assume throughout this work that the R2 utilities of interest are well-defined in the sense that integrals over the scalarisation parameters exist and are finite. For convenience, we will denote the maximum feasible scalarised value and utility by $S_{\boldsymbol{\theta}}(f(\mathbb{X})) := \max_{\mathbf{x} \in \mathbb{X}} s_{\boldsymbol{\theta}}(f(\mathbf{x}))$ and $U(f(\mathbb{X})) := \mathbb{E}_{p(\boldsymbol{\theta})}[\max_{\mathbf{x} \in \mathbb{X}}s_{\boldsymbol{\theta}}(f(\mathbf{x}))]$, respectively.

\paragraph{R2 metric} In the original work by Hansen and Jaszkiewicz \cite{hansen1998trimm}, the R2 metric, $I^{\text{R2}}: \mathbb{B}(\mathbb{R}^M) \times \mathbb{B}(\mathbb{R}^M) \rightarrow \mathbb{R}$, was defined as the difference in R2 utility between a set of objectives $Y \in \mathbb{B}(\mathbb{R}^M)$ and a user-specified reference set $Y_{R} \in \mathbb{B}(\mathbb{R}^M)$, that is
\begin{equation}
	I^{\text{R2}}(Y, Y_{R}) = U(Y_{R}) - U(Y).
	\label{eqn:r2_metric}
\end{equation}
In other words, we can interpret the R2 metric as a notion of regret for the R2 utility. Consequently, all the results that we develop in what follows for the R2 utility can also be applied for the R2 metric. Specifically, in \cref{sec:optimisation_problem}, we study the problem of maximising an R2 utility, but equivalently we could have considered the problem of minimising an R2 metric because the initial term $U(Y_{R}) \in \mathbb{R}$ is simply a constant. 

\paragraph{Qualitative aspects} The R2 utilities turn out to be a very interpretable class of multi-objective performance metrics. In particular, we can relate each component of these utility functions with the three main qualitative aspects that we described earlier in \cref{sec:utility_perspective}.

\begin{enumerate}
	\item \textbf{Convergence.} This aspect is determined by the choice of scalarisation function, which assesses the quality of each scalarised optimisation problem. To see this clearly, we can rewrite the difference between the total possible utility and the attained utility by the following expectation:
	\begin{equation*}
		U(f(\mathbb{X})) - U(f(X)) = \mathbb{E}_{p(\boldsymbol{\theta})}[
		S_{\boldsymbol{\theta}}(f(\mathbb{X})) -  S_{\boldsymbol{\theta}}(f(X))
		].
	\end{equation*}
	The difference within this expectation denotes the regret for each scalarised problem. Notably, a small regret indicates a high convergence for the scalarised problem. As the maximum utility, $U(f(\mathbb{X})) \in \mathbb{R}$, is simply a constant, we can interpret large utility values as an indication of high convergence.
	\item \textbf{Diversity.} This aspect is determined by the probability distribution, which controls the relative importance placed on each scalarised optimisation problem. As a result, we can potentially capture and assess different notions of spread and distribution by altering this density appropriately. 
	\item \textbf{Cardinality.} This aspect is determined by the $\max$ operation which ensures the utility function is monotone and therefore the addition of more points will never decrease the utility. Note that the utility function would also be monotone if we replaced the max operation with the sum operation. The main shortcoming of the summing approach is that it favours a larger set containing many poor performing points over a smaller set of high performing points.
\end{enumerate}

\paragraph{Useful properties} Another useful result about the R2 utility functions is that they satisfy both the monotone property and the diminishing returns property---we state this result in \cref{prop:monotone_submodular} and prove it in \cref{app:proofs:prop:monotone_submodular}. Intuitively, both of these properties seem sensible for a performance metric. The monotone property suggests that the quality of a set of vectors should not deteriorate as we add more points, whilst the submodularity property suggests that the improvement in adding more points should diminish with respect to the size of the set. 

Another desirable property for an R2 utility function is Pareto compliancy. This property turns out to be dependent only on the choice of scalarisation functions. In \cref{prop:weakly_pareto_compliant}, we give a sufficient condition for weak Pareto compliancy, which we prove in \cref{app:proofs:prop:weakly_pareto_compliant}. Specifically, if the scalarisation functions preserves the weak Pareto partial ordering, then so does the corresponding R2 utility. Establishing a sufficient condition for strict Pareto compliancy turns out to be non-trivial because we would have to make additional measurability assumptions. More precisely, we would have to show that the set of scalarisation parameters $\{\boldsymbol{\theta} \in \Theta: S_{\boldsymbol{\theta}}(A) > S_{\boldsymbol{\theta}}(B)\}$ has non-zero measure for any sets $A, B \in \mathbb{B}(\mathbb{R}^M)$ such that $A \succ B$.

\begin{proposition}
	All R2 utilities are monotone and submodular.
	\label[proposition]{prop:monotone_submodular}
\end{proposition}

\begin{proposition}
	If the scalarisation functions $s_{\boldsymbol{\theta}}: \mathbb{R}^M \rightarrow \mathbb{R}$ are monotonically increasing for all $\boldsymbol{\theta} \in \Theta$, then the resulting R2 utility is weakly Pareto compliant for any choice of probability distribution $p(\boldsymbol{\theta}) \geq 0$.
	\label[proposition]{prop:weakly_pareto_compliant}
\end{proposition}

\begin{remark}
	[Single-objective setting] The classical objective for single-objective optimisation is also an R2 utility, namely, $U[Y]=\max_{y \in Y} y$. As a result, all of the properties that we have shown here also applies to the standard single-objective criterion.
\end{remark}
\subsection{Special cases}
\label{sec:special_cases}
\begin{figure}
	\includegraphics[width=1\linewidth]{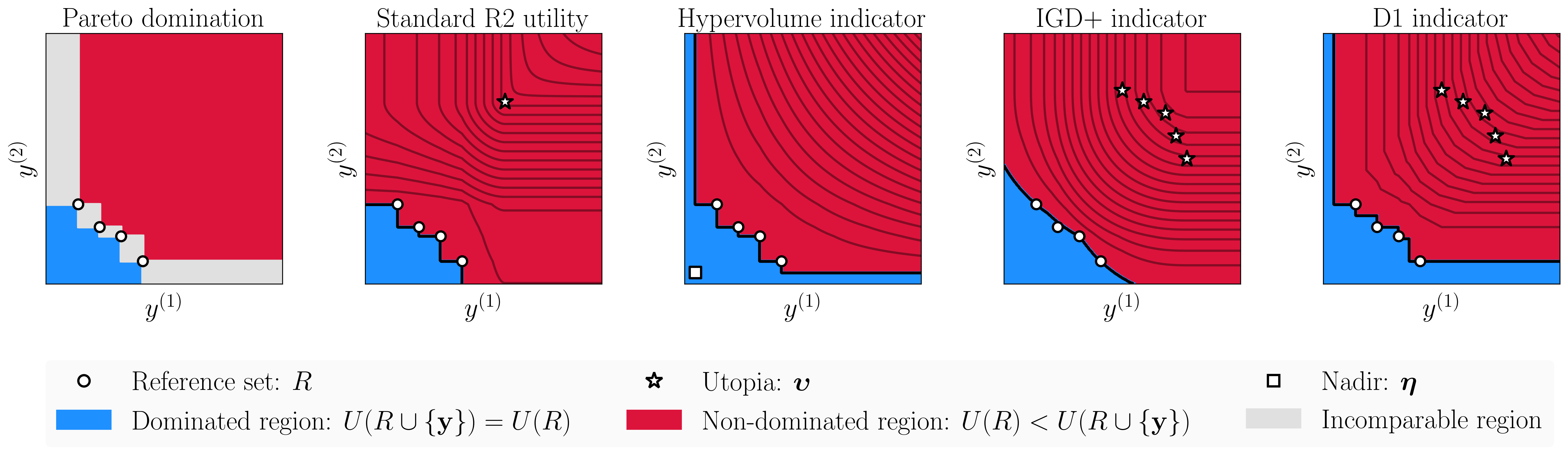}
	\caption{A comparison of the improvement regions based on the standard Pareto partial ordering over sets and some different R2 utilities in the two-objective setting. For the IGD+ utility, we set $p=2$ and $q=1$, whilst for the D1 utility we set $\mathbf{w} = (1/M, \dots, 1/M)$.}
	\label{fig:utility_domination}
\end{figure}

We now highlight the fact that many of the most widely-used performance metrics in multi-objective optimisation can be viewed as special cases of R2 utilities. In particular, from a recent survey \cite{riquelme20152laccc} of the most used performance metrics, we see that the standard R2 metric, hypervolume indicator, inverted generational distance and D1 indicator were among the top in the list. All of these utilities turn out to be special cases of R2 utilities. In \cref{fig:utility_domination}, we illustrate the contours of these R2 utilities for a simple two-dimensional example.
\subsubsection{Standard R2 metric}
In the original work by Hansen and Jaszkiewicz \cite{hansen1998trimm}, the R2 metric \eqref{eqn:r2_metric} was defined explicitly using the average weighted Chebyshev distance \eqref{eqn:chebyshev_scalarisation}
\begin{equation}
	U^{\text{R2}}(Y) = \mathbb{E}_{\mathbf{w} \sim \text{Uniform}(\Delta^{M-1})}
	\biggl[\max_{\mathbf{y} \in Y} s^{\text{Chb}}_{(\boldsymbol{\upsilon}, \mathbf{w})}(\mathbf{y}) \biggr]
	\label{eqn:standard_r2_utility}
\end{equation}
for some user-specified reference point $\boldsymbol{\upsilon} \in \mathbb{R}^M$. The original motivation behind the use of the Chebyshev scalarisation function came from the well-known result that we described in \cref{rem:chebyshev}. The use of a uniform distribution over the simplex weights was mostly done for convenience. As noted by other researchers \cite{wagner2013emo}, potentially we could have used some other weight distributions if we were interested in assessing specific parts of the Pareto front.

\paragraph{Augmented variant} A limitation of the Chebyshev scalarisation function is that it also targets weakly Pareto optimal solutions, whereas in many cases we are only interested in the strictly Pareto optimal solutions. One possible strategy to avoid these weakly Pareto optimal solution is to incorporate an $L^1$-penalty into the Chebyshev scalarisation function \cite[Part 2, Section 3.4.5]{miettinen1998}. For instance, Zitzler et al. \cite{zitzler2008moiaea} considered using the augmented Chebyshev function \eqref{eqn:augmented_chebyshev_scalarisation} instead.

\paragraph{Weak Pareto compliancy} As a result of \cref{prop:weakly_pareto_compliant}, we can conclude that both of these variants of the standard R2 utility are weakly Pareto compliant because their corresponding scalarisation functions are monotonically increasing.

\subsubsection{Hypervolume indicator} 
\label{sec:hypervolume_indicator}
One of the most popular performance criterion in multi-objective optimisation is the hypervolume indicator \cite{zitzler1998ppsn}. The hypervolume indicator, $I^{\text{HV}}: \mathbb{B}(\mathbb{R}^M) \times \mathbb{R}^M \rightarrow \mathbb{R}$, is a function that computes the volume of the dominated region between a set of objectives $Y \in \mathbb{B}(\mathbb{R}^M)$ and a reference point $\boldsymbol{\eta} \in \mathbb{R}^M$, that is
\begin{equation}
	I^{\text{HV}}(Y, \boldsymbol{\eta}) = \int_{\mathbb{R}^M} \mathbbm{1}[\mathbf{z} \in \cup_{\mathbf{y} \in Y} [\boldsymbol{\eta}, \mathbf{y}]] d\mathbf{z},
	\label{eqn:hypervolume_indicator}
\end{equation}
where $\mathbbm{1}$ is the indicator function.

\paragraph{Strict Pareto compliancy} A larger hypervolume is preferred over a smaller one. The hypervolume indicator is known to be a strictly Pareto compliant utility \cite{zitzler2003itec}, when we restrict it to sets which lie in the space that strictly Pareto dominates the reference point: $Y \in \mathbb{B}(\mathbb{D}_{\succ}(\{\boldsymbol{\eta}\})) \subset \mathbb{B}(\mathbb{R}^M)$. 

\paragraph{Hypervolume scalarisation function} Early work on the R2 indicator \cite{brockhoff2012p1acgec} identified many similarities between the behaviour of the standard R2 utility and that of the hypervolume indicator. As the hypervolume indicator grew in its popularity, a lot of work has been focussed on identifying faster and more efficient strategies to compute this indicator. After many different efforts over the last decade \cite{ishibuchi20092icec, ma2018itec, shang2018pgecc, deng2019itec, zhang2020icml}, it was slowly realised that we can write the hypervolume indicator as an R2 utility. This result is summarised in \cref{prop:hypervolume_indicator} and the proof is presented in the references \cite[Section 3.2]{shang2018pgecc}, \cite[Section 2]{deng2019itec}, and \cite[Lemma 5]{zhang2020icml}. Intuitively, the main idea of the proof is to rewrite the hypervolume integral using polar coordinates. Under this formulation, the hypervolume scalarisation function \eqref{eqn:hypervolume_scalarisation} works by computing the volume contribution over each angle and the resulting R2 utility computes the total integral of all of these contributions.

\begin{proposition}
	[Hypervolume utility] The hypervolume indicator \cref{eqn:hypervolume_indicator} can be written as an R2 utility
	\begin{equation*}
		U^{\textnormal{HV}}(Y) 
		:= I^{\textnormal{HV}}(Y, \boldsymbol{\eta}) 
		= \mathbb{E}_{\boldsymbol{\lambda} \sim \textnormal{Uniform}(\mathcal{S}_+^{M-1})}
		\biggr[\max_{\mathbf{y} \in Y} s^{\textnormal{HV}}_{(\boldsymbol{\eta}, \boldsymbol{\lambda})}(\mathbf{y}) \biggl]
	\end{equation*}
	for any $Y \in \mathbb{B}(\mathbb{R}^M)$ and $\boldsymbol{\eta} \in \mathbb{R}^M$, where the hypervolume scalarisation function is defined as the following transformation of the Chebyshev scalarisation function \eqref{eqn:chebyshev_scalarisation},
	\begin{equation}
		s^{\textnormal{HV}}_{(\boldsymbol{\eta}, \boldsymbol{\lambda})}(\mathbf{y})
		= \frac{\pi^{M/2}}{2^M \Gamma(M/2 + 1)}\min_{m=1,\dots,M} \biggl( \frac{\max(0, y^{(m)} - \eta^{(m)})}{\lambda^{(m)}} \biggr)^M,
		\label{eqn:hypervolume_scalarisation}
	\end{equation}
	with $\Gamma(z) = \int_0^\infty t^{z-1} e^{-t} dt$ denoting the Gamma function and $\boldsymbol{\lambda} \in \mathcal{S}_+^{M-1} := \{\mathbf{y} \in \mathbb{R}^M_{>0}: ||\mathbf{y}||_{L^2} = 1\}$ denoting a weight parameter that is distributed uniformly over the space of positive unit vectors.
	\label{prop:hypervolume_indicator}
\end{proposition}

\paragraph{Monotone and submodular} The hypervolume indicator is known to be a monotone and submodular set function \cite{ulrich2012lio}. These two properties have been instrumental in the development of many results and algorithms concerning the hypervolume optimisation problem \cite{guerreiro2021acs}. In \cref{prop:monotone_submodular}, we generalised this result and demonstrated that these two properties actually hold for any R2 utility. As a result, many of these existing ideas can now be applied to the R2 utility optimisation problem \eqref{eqn:utility_problem}---more details are given in \cref{sec:optimisation_problem}.
\subsubsection{Inverted generational distance} 
\label{sec:inverted_generational_distances}
Another popular performance metric in multi-objective optimisation is the inverted generational distance (IGD) \cite{coello2004m2aai}. This performance metric tries to measure some notion of distance between a finite set of objectives $Y \in \mathbb{B}(\mathbb{R}^M)$ and a finite set of ideal points $\Upsilon \in \mathbb{B}(\mathbb{R}^M)$. Using the formulation by Schutze et al. \cite{schutze2012itec}, the IGD indicator, $I^{\text{IGD}_{p, q}}: \mathbb{B}(\mathbb{R}^M) \times \mathbb{B}(\mathbb{R}^M) \rightarrow \mathbb{R}$, is equal to
\begin{equation}
	I^{\text{IGD}_{p, q}}(Y, \Upsilon) = \biggl(\frac{1}{|\Upsilon|} 
	\sum_{\boldsymbol{\upsilon} \in \Upsilon} \Bigl(
	\min_{\mathbf{y} \in Y} ||\boldsymbol{\upsilon} - \mathbf{y}||_{L^p}
	\Bigr)^q 
	\biggr)^{1/q}
	\label{eqn:igd}
\end{equation}
for some norms $p, q\geq1$. A smaller IGD is preferred because it indicates that the set of objective vectors is close to the ideal reference set. The IGD is minimised for any set of vectors containing the ideal points: $Y \supseteq \Upsilon$. It is common to set the outer norm to $q=1$, which means that the resulting IGD indicator can be interpreted as the average $L^p$-norm to the set of ideal points. On the other hand, the inner norm is commonly set to $p=2$.

\paragraph{Original formulation} The traditional definition of the IGD indicator assumes that the average is computed outside of the $L^q$-norm: $|\Upsilon|^{1/q-1} I^{\text{IGD}_{p, q}}(Y, \Upsilon)$. As described by Schutze et al. \cite[Section 3]{schutze2012itec}, the additional multiplicative factor in the traditional definition can lead to some undesirable properties. Note that there is no difference between these two formulations when $q=1$.

\paragraph{Average Lp distance} By applying a transformation on the Lp scalarisation function \eqref{eqn:lp_scalarisation}, we can rewrite the IGD indicator as a transformation of an R2 utility. We summarise this result in \cref{prop:igd} and prove it in \cref{app:proofs:prop:igd}.

\begin{proposition}
	[IGD utility] The IGD indicator \eqref{eqn:igd} can be written as a transformation of an R2 utility
	\begin{equation*}
		U^{\textnormal{IGD}_{p, q}}(Y) 
		:= -(I^{\textnormal{IGD}_{p, q}}(Y, \Upsilon))^q 
		= \frac{1}{|\Upsilon|} \sum_{\boldsymbol{\upsilon} \in \Upsilon} \max_{\mathbf{y} \in Y} s^{\textnormal{IGD}_{p, q}}_{\boldsymbol{\upsilon}}(\mathbf{y})
	\end{equation*}
	for any $Y, \Upsilon \in \mathbb{B}(\mathbb{R}^M)$ and $p, q\geq1$, where $s^{\textnormal{IGD}_{p, q}}_{\boldsymbol{\upsilon}}(\mathbf{y}) = - ||\boldsymbol{\upsilon} - \mathbf{y}||_{L^p}^q$ is the IGD scalarisation function.
	\label{prop:igd}
\end{proposition}

\paragraph{Weak Pareto compliancy} Since IGD scalarisation function is not monotonically increasing over the whole objective space, this means that the IGD indicator (and utility) is not necessarily weakly Pareto compliant. To address this problem, Ishibuchi et al. \cite{ishibuchi2015emo} proposed the IGD+ indicator, which works by truncating the difference within the $L^p$-norm at zero: $s^{\text{IGD}_{p, q}+}_{\boldsymbol{\upsilon}}(\mathbf{y}) = - ||\max(\boldsymbol{\upsilon} - \mathbf{y}, \mathbf{0}_M)||_{L^p}^q$, where the maximum is computed element-wise. This modification ensures that the scalarisation function is monotonically increasing over the entire objective space and therefore the resulting indicator (and utility) is weakly Pareto compliant.
\subsubsection{D1 indicator}
\label{sec:d1_indicator}
Another popular performance metric is the D1 indicator \cite{czyzzak1998jmda}. Similar to the IGD, the goal of the D1 indicator is to measure the distance between a finite set of objectives $Y \in \mathbb{B}(\mathbb{R}^M)$ and a finite set of ideal points $\Upsilon \in \mathbb{B}(\mathbb{R}^M)$. The D1 indicator, $I^{\text{D1}}: \mathbb{B}(\mathbb{R}^M) \times \mathbb{B}(\mathbb{R}^M) \times \Delta^{M-1} \rightarrow \mathbb{R}$ is defined as the average weighted Chebyshev distance between the set of points
\begin{equation}
	I^{\text{D1}}(Y, \Upsilon, \mathbf{w}) = \frac{1}{|\Upsilon|} \sum_{\boldsymbol{\upsilon} \in \Upsilon}
	\min_{\mathbf{y} \in Y} \max_{m=1,\dots,M} w^{(m)} (\upsilon^{(m)} - y^{(m)}),
	\label{eqn:d1}
\end{equation}
where $\mathbf{w} \in \Delta^{M-1}$ denotes some weight vector. A smaller D1 indicator is preferred because it indicates that the set of objectives is close to the set of ideal points. It is common to see uniform weights being used for each objective: $\mathbf{w} = \mathbf{1}_M / M \in \Delta^{M-1}$. 

\paragraph{Average Chebyshev distance} The negative D1 indicator can be obtained as an R2 utility using the Chebyshev scalarisation function \eqref{eqn:chebyshev_scalarisation} and a uniform distribution over the set of ideal points. We summarise this result in \cref{prop:d1} and prove it in \cref{app:proofs:prop:d1}.

\begin{proposition}
	[D1 utility] The D1 indicator \eqref{eqn:d1} can be written as a transformation of an R2 utility
	\begin{equation*}
		U^{\textnormal{D1}}(Y) := -I^{\textnormal{D1}}(Y, \Upsilon, \mathbf{w}) 
		= \frac{1}{|\Upsilon|} \sum_{\boldsymbol{\upsilon} \in \Upsilon}
		\max_{\mathbf{y} \in Y} s^{\textnormal{Chb}}_{(\boldsymbol{\upsilon}, \mathbf{w})}(\mathbf{y})
	\end{equation*}
	for any $Y, \Upsilon \in \mathbb{B}(\mathbb{R}^M)$ and $\mathbf{w} \in \Delta^{M-1}$.
	\label{prop:d1}
\end{proposition}

\paragraph{Weak Pareto compliancy} The D1 indicator (and utility) is weakly Pareto compliant because the Chebyshev scalarisation function is monotonically increasing. 
\subsection{Discussion}
\label{sec:r2_discussion}
The overall goal of multi-objective optimisation is to identify an approximation to the Pareto optimal points. The scalarisation and utility perspectives of multi-objective optimisation offers a strategy to obtain this approximation by recasting the original problem as one which is more amenable to standard optimisation strategies. The R2 utilities unify both of these approaches into one in which a decision maker only needs to specify a family of scalarisation functions and a probability distribution to go along with it. There is no clear-cut answer as to what utilities, scalarisation functions and probability distribution we should use in practice. This choice naturally depends on the goals and preferences of the decision maker, which can in general be hard to elicit \cite{wang2017cis}. 

In many cases, decision makers often resort to the use of popular off-the-shelf utility functions such as those listed in \cref{sec:special_cases}. These utility functions are useful for a general decision maker, who is simply interested in obtaining some approximation to the Pareto optimal points, but they might be insufficient for a more specialised user who has specific preferences they want to encode in their approximation. For example, suppose a decision maker has a lexicographical preference on the objectives, which states that some objectives are more important to optimise than others. The previously listed performance metrics might be inadequate for this setting because they do not explicitly encode this preference into their utility values. A more suitable approach would be to design an R2 utility which does encode this preference. For instance, we could use a weight-based scalarisation function such as the Chebyshev scalarisation function \eqref{eqn:chebyshev_scalarisation} and devise a weight distribution that satisfies this lexicographical constraint.

In this paper, we have highlighted only a small selection of R2 utilities in order to showcase the usefulness of the methodology. There are numerous possibilities that we have not discussed. For instance, we could consider R2 utilities defined using reference-line based scalarisation functions \cite{das1998sjo,zhang2007itec,sato2014p2acgec,namura2017itec}. Similarly, we could consider R2 utilities based on transformations of existing scalarisation functions such as the linear or Chebyshev scalarisation function \cite{marler2004smo,ishibuchi20092icec,ma2018itec,binois2020jmlr,chugh20202icecc}.
\section{Optimisation problem}
\label{sec:optimisation_problem}
Suppose that we have identified an R2 utility which adequately reflects the decision maker's preferences. We will now consider solving the corresponding utility optimisation problem 
\begin{equation}
	\max_{X \subseteq \mathbb{X}, |X| \leq P} \mathbb{E}_{p(\boldsymbol{\theta})}\biggl[\max_{\mathbf{x} \in X} s_{\boldsymbol{\theta}}(f(\mathbf{x})) \biggr]
	\label{eqn:r2_utility_problem}
\end{equation}
for some $P > 0$. As all R2 utilities are monotone and submodular, this utility optimisation problem turns out to be an instance of a cardinality-constrained submodular optimisation problem, which is known to be NP-hard in general \cite{bach2013fiml, krause2014tpathp}. Nevertheless, we can typically solve these problems approximately using greedy algorithms. These greedy approaches satisfy the well-known performance guarantee by Nemhauser, Wolsey and Fisher \cite{nemhauser1978mp}. In \cref{sec:optimisation_strategies}, we will review how these greedy strategies can be used to approximately solve our utility optimisation problem. We will then move on to \cref{sec:bayesian_optimisation}, where we will extend the analysis of these greedy strategies to the Bayesian optimisation setting.

In order to prove some later results, we will now assume for convenience that the R2 utility function of interest is non-negative and bounded over the feasible space. Note that any R2 utility can be made non-negative over a bounded space by simply adding a large enough constant.
\begin{assumption}
	[Non-negative and bounded] Consider an objective function $f: \mathbb{X} \rightarrow \mathbb{R}^M$ and an R2 utility $U: \mathbb{B}(\mathbb{R}^M) \rightarrow \mathbb{R}$. We assume that there exists a constant $C>0$ such that $0 \leq s_{\boldsymbol{\theta}}(f(\mathbf{x})) \leq C$ for all $\mathbf{x} \in \mathbb{X}$ and $\boldsymbol{\theta} \in \Theta$.
	\label{ass:non-negative_and_bounded}
\end{assumption}
\subsection{Greedy algorithms}
\label{sec:optimisation_strategies}
To solve the R2 utility optimisation problem \eqref{eqn:r2_utility_problem} approximately, we will consider the use of greedy algorithms. These algorithms start with an empty set of inputs $X_0 = \{\}$, which is then incrementally augmented according to some greedy heuristic: $X_n = X_{n-1} \cup \{\mathbf{x}_n\}$ for $n=1, \dots, N$, where $\mathbf{x}_n \in \mathbb{X}$ is the point that is greedily selected at the $n$-th round. The final set $X_N$ is then suggested as a solution, where $N$ is the budget of function evaluations.
\subsubsection{Standard greedy approach} 
In the traditional greedy strategy, we pick the input which achieves the largest utility improvement; that is,
\begin{align}
	\begin{split}
		\mathbf{x}_{n+1} 
		&\in \argmax_{\mathbf{x} \in \mathbb{X}} (U(f(X_n \cup \{\mathbf{x}\})) - U(f(X_n))).
	\end{split}
	\label{eqn:greedy_strategy}
\end{align}
This greedy method comes with the well-known performance guarantee by Nemhauser et al. \cite{nemhauser1978mp}, which we recall in \cref{thm:greedy_guarantee}. The main consequence of this guarantee is that after $P$ rounds of the greedy algorithm, the utility obtained is guaranteed to be at least $(1-e^{-1})$ times the best possible utility for this problem. Moreover, if we are allowed to perform even more rounds, then this approximation becomes even better. For instance, after $\gamma P$ rounds, the factor is equal to $(1-e^{-\gamma})$, which becomes ever closer to one as $\gamma$ grows.
\begin{theorem}
	[Greedy guarantee] (Nemhauser, Wolsey and Fisher \cite{nemhauser1978mp}) Consider an objective function $f:\mathbb{X} \rightarrow \mathbb{R}^M$ and an R2 utility $U: \mathbb{B}(\mathbb{R}^M) \rightarrow \mathbb{R}$, satisfying \cref{ass:non-negative_and_bounded}. Let $\{X_n\}_{n\geq1}$ be the greedily selected inputs, then for all positive integers $P$ and $N$,
	\begin{equation*}
		U(f(X_N)) \geq (1 - e^{-N/P}) \max_{X \subseteq \mathbb{X}: |X| \leq P} U(f(X)).
		\label{eqn:greedy_guarantee}
	\end{equation*}
	\label{thm:greedy_guarantee}
\end{theorem}

\subsubsection{Approximate greedy approach}
The R2 utility can only be evaluated exactly in simple settings such as the case where the set of scalarisation parameters is finite. In general, we rely on approximations of the R2 utility such as the one obtained via Monte Carlo estimation, that is
\begin{equation}
	\hat{U}_J(Y) = \frac{1}{J} \sum_{j=1}^J \max_{\mathbf{y} \in Y}s_{\boldsymbol{\theta}_j}(\mathbf{y}),
	\label{eqn:utility_estimate}
\end{equation}
where $\boldsymbol{\theta}_j \sim p(\boldsymbol{\theta})$ denotes the $J$ independent samples of the scalarisation parameter. Equipped with this estimate, the approximate greedy strategy proceeds by maximising the estimate of the utility improvement
\begin{align}
	\begin{split}
		\mathbf{x}_{n+1} 
		&\in \argmax_{\mathbf{x} \in \mathbb{X}} (\hat{U}_J(f(X_n \cup \{\mathbf{x}\})) - \hat{U}_J(f(X_n))).
	\end{split}
	\label{eqn:approximate_greedy_strategy}
\end{align}
Note that the estimate of the R2 utility \eqref{eqn:utility_estimate} could be fixed throughout the duration of the optimisation procedure or could be re-estimated at every iteration. In addition, we could also let the number of samples change over time as well---for simplicity we only consider the static case where the number of samples is fixed throughout.

\paragraph{Random scalarisation} When we re-estimate the utility at each time and use only one Monte Carlo sample, $J=1$, then we recover the random scalarisation algorithm. This is a popular strategy that has been used before in the area of Bayesian optimisation \cite{knowles2006itec, paria2020uai, zhang2020icml,chowdhury2021icais} and bares some similarity with Thompson sampling strategies \cite{thompson1933b, russo2014moo}.

\paragraph{Approximate greedy guarantee} By a simple application of Hoeffding's inequality, we can devise a similar performance guarantee to the one above. The main difference is that we now have an additional additive penalty term, which arises from the error that we incur in performing Monte Carlo estimation at every round. This error naturally decreases as the number of Monte Carlo samples increases. For completeness, we state this result in \cref{thm:approximate_greedy_guarantee} and prove it in \cref{app:proofs:thm:approximate_greedy_guarantee}.
\begin{proposition}
	[Approximate greedy guarantee] Consider an objective function $f: \mathbb{X} \rightarrow \mathbb{R}^M$ and an R2 utility $U: \mathbb{B}(\mathbb{R}^M) \rightarrow \mathbb{R}$, satisfying \cref{ass:non-negative_and_bounded}. Let $\{X_n\}_{n\geq1}$ denote the set of inputs identified by the approximately greedy strategy using $J>0$ Monte Carlo samples, then for any $\delta \in (0, 1)$ and any positive integers $P$ and $N$, the following inequality holds with probability $1-\delta$:
	\begin{align*}
		U(f(X_N)) 
		&\geq (1-e^{-N/P}) \max_{X \subseteq \mathbb{X}, |X|\leq P} U(f(X)) - \epsilon(\delta, P, N, J),
		\label{eqn:approximate_greedy_guarantee}
	\end{align*}
	where the error term is given by $\epsilon(\delta, P, N, J) = \sqrt{\frac{2}{J}  \log(\frac{4N}{\delta})} \sum_{n=1}^N C_{n-1} (1- \frac{1}{P})^{N-n}$, with $C_n := \sup_{\boldsymbol{\theta} \in \Theta} (S_{\boldsymbol{\theta}}(f(\mathbb{X})) - S_{\boldsymbol{\theta}}(f(X_n)))$ denoting the maximum scalarised regret at round $n=1,\dots,N-1$ and $C_0 := C$.
	\label{thm:approximate_greedy_guarantee}
\end{proposition}

\paragraph{Number of samples} As expected, the error decreases to zero as we increase the number of Monte Carlo samples: $\epsilon(\delta, P, N, J) \rightarrow 0$ as $J\rightarrow \infty$. Note that we can bound this error by an input-independent constant $\epsilon(\delta, P, N, J) \leq CP \sqrt{2 \log(4N/\delta) / J}$. As a result, if we wanted to achieve a specific error of $\epsilon > 0$ with probability $1-\delta$, then we could set the number of Monte Carlo samples to be $J = \lceil 2C^2 P^2 \log(4N/\delta) / \epsilon^2 \rceil$. This number is arguably much larger than what is required in practice and could be easily improved by a more careful application of some stronger concentration inequalities. Nevertheless, this result gives a useful indication of what factors affect the performance guarantee. For instance, we see that the number of objectives $M$ only enters into the calculation via the upper bound $C$, which is determined by both the objective function and the R2 utility.

\paragraph{Noisy submodular optimisation} The idea of solving a submodular optimisation problem when there is uncertainty in the function evaluations has been addressed before in the area of submodular optimisation under noise. For example, Kempe, Kleinbery and Tardos \cite{kempe2003pnasickddm} considered an influence optimisation problem in a social network, whilst Singla, Tschiatschek, and Krause \cite{singla2016ptacai} considered some applications for crowdsourced image collection. The fact that we incur an additive error in the performance bound is known to be a direct consequence of the uncertainty in the function evaluation. In our case, this error arises from our Monte Carlo estimate of the R2 utility, which can be easily reduced at the expense of more computational power. 

We now illustrate a simple example of this approximate greedy strategy on the hypersphere problem in \cref{eg:hypersphere}. In particular, we illustrate empirically how the observed performance changes in this problem when we vary the number of Monte Carlo samples and the number of objectives.

\begin{figure}
	\includegraphics[width=0.85\linewidth]{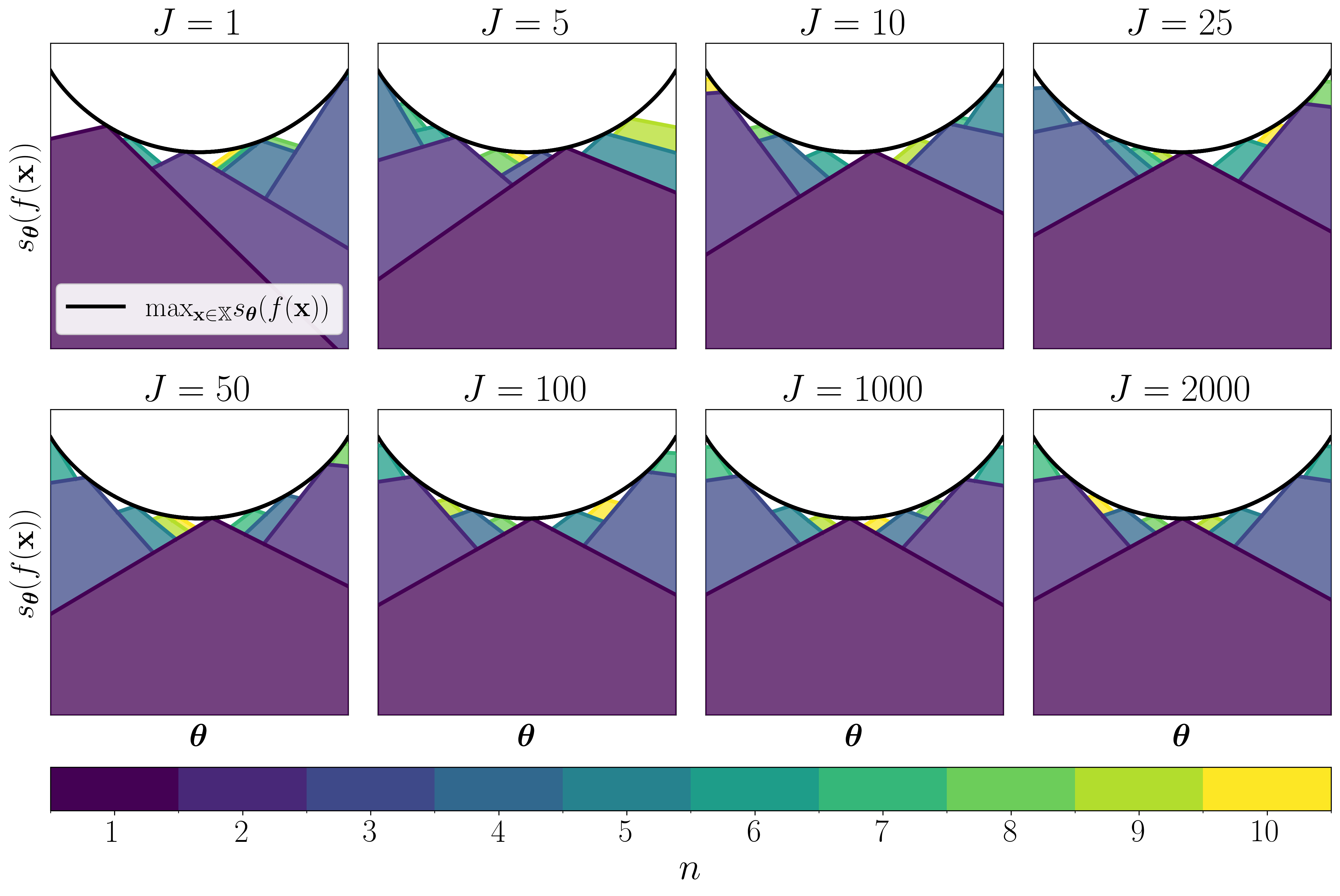}
	\centering
	\caption{A visual comparison of the approximate greedy strategy with different number of Monte Carlo samples on the hypersphere problem (\cref{eg:hypersphere}).}
	\label{fig:approximate_greedy_mc}
\end{figure}

\begin{figure}
	\includegraphics[width=1\linewidth]{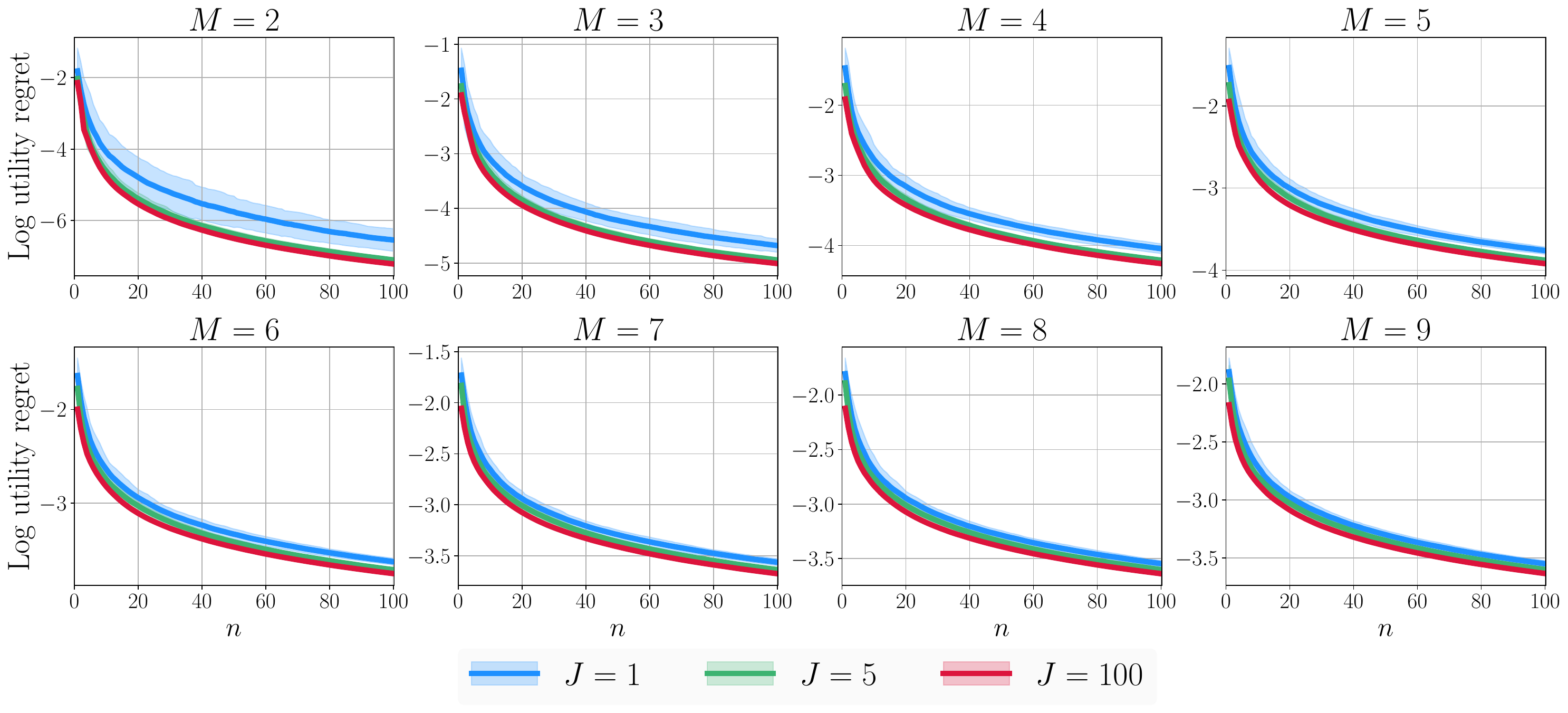}
	\centering
	\caption{A performance comparison of the approximate greedy strategy when varying the number of Monte Carlo samples and the number of objectives on the hypersphere problem (\cref{eg:hypersphere}). The mean and two standard deviations of the log utility regret, $\log(U(f(\mathbb{X})) - U(f(X_N)))$, obtained from 100 independent runs, are plotted. In each experiment, the Pareto front and the standard R2 utility were approximated using $10^5$ points and samples.} 
	\label{fig:utility_optimisation_m}
\end{figure}

\begin{example}
	[Hypersphere] Consider the R2 utility optimisation problem \eqref{eqn:r2_utility_problem}, where the feasible objective space $f(\mathbb{X})$ is the $M$-dimensional hypersphere and the performance metric $U$ is the standard R2 utility \eqref{eqn:standard_r2_utility} with the reference point set to $\boldsymbol{\upsilon}=(1,\dots,1) \in \mathbb{R}^M$. For this problem, the Pareto front of interest is just the surface of the hypersphere lying in the non-negative orthant: $\mathbb{Y}^* = \mathcal{S}^{M-1}_+$, which we defined in \cref{sec:hypervolume_indicator}. 
	
	In \cref{fig:approximate_greedy_mc}, we illustrate one run of the approximate greedy strategy on this problem for the $(M=2)$-dimensional setting with different numbers of Monte Carlo samples $J$. We plot the result in the scalarisation parameter space. To elaborate, at each time $n$, we picked a solution $\mathbf{x}_n \in \mathbb{X}$ according to the approximate greedy heuristic \eqref{eqn:approximate_greedy_strategy}. This solution is associated with a surface in the scalarisation parameter space: $\{s_{(\boldsymbol{\upsilon}, (w, 1-w))}^{\text{Chb}}(f(\mathbf{x}_n)): w \in [0, 1]\}$. In this illustrative example, the scalarisation parameter can be parametrised using a one-dimensional weight, $\boldsymbol{\theta} = w$, which lies in the unit interval $\Theta = [0, 1]$. Geometrically, we can interpret the R2 utility of this solution, $U(\{f(\mathbf{x}_n)\}) \in \mathbb{R}$, as the area underneath this curve. Moreover, we can interpret the utility of a collection of these solutions as the area underneath the ``envelope'' of these curves. In the figure, we coloured the gain in utility at each time based on the iteration number $n$. We see clearly that the gain in utility is much larger at earlier rounds, which is a direct consequence of the diminishing returns property. In addition, we see how the number of Monte Carlo samples effects the points that are selected. Specifically, when the number of samples is small, then there is a large variance in the points that are selected---however, this variance however decreases when the number of samples increases.
	
	In \cref{fig:utility_optimisation_m}, we illustrate the performance of the approximate greedy strategy when we vary the number of objectives as well as the number of Monte Carlo samples.  As expected from the diminishing return property, we see that the improvement in utility drops quickly as the number of iterations increases. Moreover, we see a small improvement in performance when we increase the number of Monte Carlo samples. 
	\label{eg:hypersphere}
\end{example}
\subsection{Bayesian optimisation}
\label{sec:bayesian_optimisation}
In the previous section, we discussed how greedy strategies can be used to solve the R2 utility optimisation problem \eqref{eqn:r2_utility_problem}. All of these strategies rely on solving a sequence of optimisation problems that depend on exact evaluations of the underlying vector-valued objective function. In this section, we extend these ideas to the black-box optimisation setting, where the function evaluations are assumed to be expensive and potentially subject to noise. We focus our study on the popular Bayesian optimisation algorithm, which relies on the use of surrogate models in order to decide which points to evaluate. We summarise the main steps of this optimisation procedure in what follows and provide the pseudo-code in \cref{alg:bayesian_optimisation}---for a more in-depth overview, consult the recent review by Garnett \cite{garnett2023}.

\paragraph{Modelling} At time $n$, we have a data set of inputs and outputs, $\mathcal{D}_n = \mathcal{D}_0 \cup \{(\mathbf{x}_t, \mathbf{y}_t)\}_{t=1,\dots,n}$. Given this data set, the Bayesian optimisation routine proceeds by building a probabilistic surrogate model of the objective function. At a high level, this model can be constructed by adopting either a Bayesian perspective or a frequentist perspective.

\begin{enumerate}[label=(i)]
	\item In the Bayesian paradigm, we assume a prior distribution on the objective function $p(f|\mathcal{D}_0)$ and a likelihood $p(\mathbf{y}| \mathbf{x}, f)$ on the observations. Then by standard conditioning we can obtain a posterior distribution over the function $p(f|\mathcal{D}_n) \propto p(f|\mathcal{D}_0) \prod_{t=1}^n p(\mathbf{y}_t|\mathbf{x}_t, f)$.
	\item In the frequentist paradigm, we directly construct an estimate of the function $\hat{f}: \mathbb{X} \rightarrow \mathbb{R}^M$. This is often done by minimising some empirical loss function which depends on the observed data points.
\end{enumerate}
To keep the discussion general, we denote the probabilistic surrogate model, at time $n$, under both settings, by $p(\hat{f}|\mathcal{D}_n)$. In the Bayesian setting, this can be interpreted as a posterior distribution of the model at time $n$, whilst in the frequentist setting, this can be interpreted as a point mass on our function estimate at time $n$. 

\paragraph{Acquisition} Equipped with a surrogate model, the Bayesian optimisation procedure then proceeds to select a new input to query by maximising some acquisition function $\alpha: \mathbb{X} \rightarrow \mathbb{R}$,
\begin{equation*}
	\mathbf{x}_{n+1} \in \argmax_{\mathbf{x} \in \mathbb{X}} \alpha(\mathbf{x} | \mathcal{D}_n).
\end{equation*}
This acquisition function is designed to strike a balance between exploring new parts of input space and exploiting the regions of the input space that are known to perform well. In multi-objective Bayesian optimisation, practitioners often rely on the use of scalarisation functions or utility functions in order to navigate this trade-off. In this paper, we show how it is possible to unify these existing approaches by considering a general class of acquisition functions based on the R2 utility.

\begin{algorithm}
	\DontPrintSemicolon
	\SetKwFunction{FMain}{\normalfont \textsc{BayesianOptimisation}}
	\SetKwProg{Fn}{Function}{:}{}
	\nonl \Fn{\FMain{$\mathcal{D}_0$, $N$, $p(\hat{f}| \mathcal{D}_0)$, $\alpha$}}{
		\nonl \textcolor{black}{\nonl \textit{\small // Initial data set $\mathcal{D}_0$.}}

		\nonl \textcolor{black}{\nonl \textit{\small // Number of function evaluations $N$.}}

		\nonl \textcolor{black}{\nonl \textit{\small // Initial model $p(\hat{f}|\mathcal{D}_0)$.}}

		\nonl \textcolor{black}{\nonl \textit{\small // The acquisition function $\alpha: \mathbb{X} \rightarrow \mathbb{R}$.}}
		
		\For{$n=0,\dots,N-1$}{
			Acquire the new input: $\mathbf{x}_{n+1} \in \argmax_{\mathbf{x} \in \mathbb{X}} \alpha(\mathbf{x}|\mathcal{D}_n)$.

			Evaluate the input: $\mathbf{y}_{n+1} \approx f(\mathbf{x}_{n+1})$.

			Update the data set $\mathcal{D}_{n+1} = \mathcal{D}_n \cup \{(\mathbf{x}_{n+1}, \mathbf{y}_{n+1})\}$.

			Update the model: $p(\hat{f}|\mathcal{D}_{n+1})$.
		}
		\Return The final data set $\mathcal{D}_N$.
	}
	\caption{}
	\label{alg:bayesian_optimisation}
\end{algorithm}
\subsubsection{Adjusted expected utility improvement}
We now introduce a general family of greedy acquisition functions called the adjusted expected utility improvement (AEUI). An acquisition function belongs to this family if it can be written in the form
\begin{equation}
	\alpha^{\text{AEUI}}(\mathbf{x}| \mathcal{D}_n) = \alpha^{\text{EUI}}(\mathbf{x}| \mathcal{D}_n) + \mathcal{A}_\delta(\mathbf{x}| \mathcal{D}_n),
	\label{eqn:aeui}
\end{equation}
where $\alpha^{\text{EUI}}: \mathbb{X} \rightarrow \mathbb{R}$ denotes the expected utility improvement (EUI) acquisition function,
\begin{align}
	\begin{split}
		\alpha^{\text{EUI}}(\mathbf{x}| \mathcal{D}_n) 
		&= \mathbb{E}_{p(\hat{f}|\mathcal{D}_n)} [ U(\hat{f}(X_n \cup \{\mathbf{x}\})) - U(\hat{f}(X_n)) ],
	\end{split}
	\label{eqn:eui}
\end{align}
and $\mathcal{A}_\delta: \mathbb{X} \rightarrow \mathbb{R}_{\geq 0}$ denotes the adjustment function, which is a non-negative function that depends on some parameter $\delta \in (0, 1)$.

\paragraph{Expected utility improvement} The EUI acquisition function, can be interpreted as the natural counterpart to the standard greedy heuristic in \eqref{eqn:greedy_strategy}. Instead of computing the greatest utility improvement with respect to the true objective function, we compute the greatest EUI with respect to the surrogate model. Notably, this family of acquisition functions has appeared numerous times before in the literature---see the recent review by Zhan and Xing \cite{zhan2020jgo}.

\paragraph{Adjustment function} $\mathcal{A}_\delta$ can be interpreted as a function that accounts for the error from estimating the objective function using the surrogate model. Specifically, we will assume that the adjustment function is an upper bound on the absolute error between the expected utility and the actual utility---we formalise this statement below in \cref{ass:concentration}.

\begin{assumption}
	[Utility estimate concentration] Consider an objective function $f: \mathbb{X} \rightarrow \mathbb{R}^M$, an R2 utility $U: \mathbb{B}(\mathbb{R}^M) \rightarrow \mathbb{R}$, a surrogate model $p(\hat{f}|\mathcal{D}_n)$ and a family of adjustment functions $\{\mathcal{A}_{\delta}: \mathbb{X} \rightarrow \mathbb{R}_{\geq 0}: \delta \in (0, 1)\}$. For any $\delta \in (0, 1)$, we assume that the following inequality holds with probability $1-\delta$:
	\begin{equation}
		\bigl|U(f(X_n \cup \{\mathbf{x}\})) - \mathbb{E}_{p(\hat{f}|\mathcal{D}_n)}[U(\hat{f}(X_n \cup \{\mathbf{x}\}))] \bigr| 
		\leq \mathcal{A}_\delta(\mathbf{x}| \mathcal{D}_n)
		\label{eqn:concentration_inequality}
	\end{equation} 
	for all inputs $\mathbf{x} \in \mathbb{X}$.
	\label{ass:concentration}
\end{assumption}

\begin{remark}
	Note that the family of EUI acquisition functions is a subset of the family of AEUI acquisition functions. More explicitly, to recover the EUI acquisition policy, we can simply set the adjustment function to be any constant function that satisfies \cref{ass:concentration}.
	\label{rem:eui}
\end{remark}

\begin{remark}
	[Batch optimisation] All of the results that we develop in this section can easily be extended to the batch optimisation setting, where we pick multiple points at a time rather than just one. More precisely, for the batch setting, we would consider the batch AEUI acquisition function
	\begin{equation*}
		\alpha^{\textnormal{AEUI}}(X| \mathcal{D}_n) 
		= \mathbb{E}_{p(\hat{f}|\mathcal{D}_n)} [ U(\hat{f}(X_n \cup X)) - U(\hat{f}(X_n)) ] + \mathcal{A}_\delta(X| \mathcal{D}_n),
	\end{equation*}
	where $X \in \mathbb{X}^q$ is a batch of $q > 0$ inputs and $\mathcal{A}_\delta(X| \mathcal{D}_n)$ is the adjustment term that satisfies the batch analogue of the concentration inequality described in \cref{ass:concentration}. For ease of exposition, we have intentionally focussed our discussion on only the sequential setting.
\end{remark}
\subsubsection{Computing the acquisition function}
To compute the AEUI acquisition function, we need to be able to evaluate both the EUI acquisition function and the adjustment function.

\paragraph{Estimating the expected utility improvement} In general, the EUI acquisition function is not a tractable quantity because we might not be able to evaluate the utility exactly or be able to compute the expectation over the surrogate model. Nevertheless, if we can generate samples of the scalarisation parameters and surrogate model, then we can compute a Monte Carlo estimate
\begin{align}
	\begin{split}
		\hat{\alpha}^{\text{EUI}}(\mathbf{x}| \mathcal{D}_n) 
		&= \frac{1}{H} \sum_{h=1}^H \hat{U}_J(\hat{f}_h(X_n \cup \{\mathbf{x}\})) - \hat{U}_J(\hat{f}_h(X_n)),
	\end{split}
	\label{eqn:approximate_eui}
\end{align} 
where $\hat{f}_h(\cdot) \sim p(\hat{f}|\mathcal{D}_n)$ denotes the independent samples from the model, whilst $\hat{U}_J$ denotes the Monte Carlo estimate of the utility described earlier in \eqref{eqn:utility_estimate}. In order to develop the performance bound in the next section, we will assume that the estimates for the utility based on the surrogate model are non-negative and bounded by the same constant, $C>0$, which bounds the actual utility values in \cref{ass:non-negative_and_bounded}.
\begin{assumption}
	[Non-negative and bounded utility estimates] Consider an R2 utility and a surrogate model $p(\hat{f}|\mathcal{D}_n)$. We assume that there exists a constant $C>0$ such that $0 \leq s_{\boldsymbol{\theta}}(\hat{f}(\mathbf{x})) \leq C$ for all $\mathbf{x} \in \mathbb{X}$ and $\boldsymbol{\theta} \in \Theta$, almost surely over all possible model samples $\hat{f}(\cdot) \sim p(\hat{f}|\mathcal{D}_n)$.
	\label{ass:non-negative_and_bounded_estimate}
\end{assumption}

\begin{remark}
	[Sampling points instead of paths] Sampling the whole path $\hat{f}_h$ is expensive and unnecessary in practice. We can simplify the estimate of the EUI acquisition function by marginalising out the dependency on the points that are not considered in the integral. As a result, we only need to sample the function at the considered locations, $\hat{f}_h(X_n \cup \{\mathbf{x}\}) \sim p(\hat{f}(X_n \cup \{\mathbf{x}\})|\mathcal{D}_n)$, which is considerably cheaper. 
	\label[remark]{rem:sampling_points}
\end{remark}

\paragraph{Computing the adjustment function} The adjustment term quantifies the error in estimating the utility of the objective function. If we are only interested in employing the EUI acquisition function, then we do not have to compute this quantity because it would only contribute an input-independent constant (\cref{rem:eui}). On the other hand, if we wanted to establish an input-dependent adjustment term, then we would need to make some additional assumptions on both the utility function and the model. For example, if we assumed that the scalarisation functions are Lipschitz in the $L^2$-norm, then we can set the adjustment function to be proportional to any upper bound on the expected $L^2$-error. To see this, suppose there exists a constant $L > 0$ such that $|s_{\boldsymbol{\theta}}(\mathbf{y}) - s_{\boldsymbol{\theta}}(\mathbf{y}')| \leq L ||\mathbf{y} - \mathbf{y}'||_{L^2}$ for any $\boldsymbol{\theta} \in \Theta$ and $\mathbf{y}, \mathbf{y}' \in \mathbb{R}^M$, then by the triangle inequality we obtain the inequality
\begin{align*}
	\bigl|U(f(X)) - \mathbb{E}_{p(\hat{f}|\mathcal{D}_n)}[U(\hat{f}(X))] \bigr|
	\leq L \mathbb{E}_{p(\hat{f}|\mathcal{D}_n)} \biggl[ \max_{\mathbf{x}' \in X} ||f(\mathbf{x}') - \hat{f}(\mathbf{x}') ||_{L^2}\biggr],
\end{align*}
which holds for any finite set of inputs $X \subseteq \mathbb{X}$. Therefore, if we can bound the expected $L^2$-error of the function estimate at $X = X_n \cup \{\mathbf{x}\}$, then we have a candidate for the adjustment function. For instance, if we adopt a standard Gaussian process \cite{rasmussen2006} regression set-up, similar to the one described by Chowdhury and Gopalan \cite{chowdhury2021icais}, then we can set the adjustment function to
\begin{equation}
	\mathcal{A}_\delta(\mathbf{x}|\mathcal{D}_n) = a_{\delta} + b_{\delta} \text{trace}(\boldsymbol{\Sigma}_n(\mathbf{x}, \mathbf{x}))^{1/2},
	\label{eqn:gp_adjustment}
\end{equation}
where $\boldsymbol{\Sigma}_n$ denotes the covariance matrix of the Gaussian process model at time $n$, whilst the variables $a_{\delta}, b_{\delta}$ are simply some model-dependent constants which are independent of the input.
\subsubsection{Theoretical guarantee}
We will now extend the greedy performance guarantees described before in \cref{sec:optimisation_strategies} for the Bayesian optimisation setting. Specifically, our main result in \cref{thm:aeui_guarantee}, gives a performance bound for the general family of approximate AEUI acquisition functions
\begin{align}
	\begin{split}
		\hat{\alpha}^{\text{AEUI}}(\mathbf{x}| \mathcal{D}_n) 
		&= \hat{\alpha}^{\text{EUI}}(\mathbf{x}| \mathcal{D}_n) + \mathcal{A}_\delta(\mathbf{x}| \mathcal{D}_n),
	\end{split}
	\label{eqn:approximate_aeui}
\end{align} 
where $\hat{\alpha}^{\text{EUI}}: \mathbb{X} \rightarrow \mathbb{R}$ is the approximate EUI acquisition function in \eqref{eqn:approximate_eui}. The proof of this result is presented in \cref{app:proofs:thm:aeui_guarantee}. 

\begin{proposition}
	[Approximate AEUI guarantee] Consider an objective function $f:\mathbb{X} \rightarrow \mathbb{R}^M$ and an R2 utility $U: \mathbb{B}(\mathbb{R}^M) \rightarrow \mathbb{R}$ such that \cref{ass:non-negative_and_bounded,ass:concentration,ass:non-negative_and_bounded_estimate} holds. Let $\{X_n\}_{n\geq 1}$ be the inputs selected using the corresponding approximate AEUI acquisition function, then for any $\delta \in (0, 1)$ and any positive integers $P$ and $N$, the following inequality holds with probability $1-\delta$:
	\begin{align}
		U(f(X_N)) 
		&\geq (1-e^{-N/P}) \max_{X \subseteq \mathbb{X}, |X|\leq P} U(f(X)) - \sum_{i=1}^2 \epsilon_i(\delta, P, N, J, H),
		\label{eqn:aeui_guarantee}
	\end{align}
	where the error is a sum of the Monte Carlo error and the model error,
	\begin{align*}
		\epsilon_1(\delta, P, N, J, H)& =  \sqrt{\frac{2}{\min(J, H)} \log\biggl(\frac{12 N}{\delta}\biggr)} \sum_{n=1}^N C_{n-1} \biggl(1-\frac{1}{P}\biggr)^{N-n},
		\\
		\epsilon_2(\delta, P, N, J, H) &= 2 \sum_{n=1}^N \mathcal{A}_{\delta/(6N)}(\mathbf{x}_n| \mathcal{D}_{n-1}) e^{-(N-n)/P},
	\end{align*}
	respectively, where $C_n := \sup_{\boldsymbol{\theta} \in \Theta}\mathbb{E}_{p(\hat{f}|\mathcal{D}_n)}[S_{\boldsymbol{\theta}}(\hat{f}(\mathbb{X})) - S_{\boldsymbol{\theta}}(\hat{f}(X_n))]$ denotes the expected maximal scalarised regret for $n=1,\dots,N-1$ and $C_0 := C$.
	\label{thm:aeui_guarantee}
\end{proposition}

\paragraph{Monte Carlo error} As reflected in the performance bound, we incur an additive error from estimating the EUI acquisition function using Monte Carlo. This error naturally decreases when we increase the number of samples. In particular, if we can compute the utility exactly, then this is the equivalent to setting $J \rightarrow \infty$. Similarly, if we can compute the expectation exactly, then this is the same as setting $H \rightarrow \infty$. Notably, if we can compute both the utility and the expectation exactly, then this error term becomes zero: $\epsilon_1(\delta, P, N, J, H) = 0$. In this case, the resulting performance bound would reflect the one obtained by employing the exact AEUI acquisition function \eqref{eqn:aeui}.

\paragraph{Model error} The second error term in the performance bound reflects the error incurred in approximating the objective function using the surrogate model at every round. Note that the errors from the earlier rounds $n$ are down-weighted exponentially. Therefore, if the model becomes more representative of the true objective function, as we acquire more data, then this error tends towards zero: $\epsilon_2(\delta, P, N, J, H) \rightarrow 0$ as $N\rightarrow \infty$. This result is expected because there should be no downside in using the surrogate model when it adequately describes the underlying objective function.

\begin{figure}
	\includegraphics[width=1\linewidth]{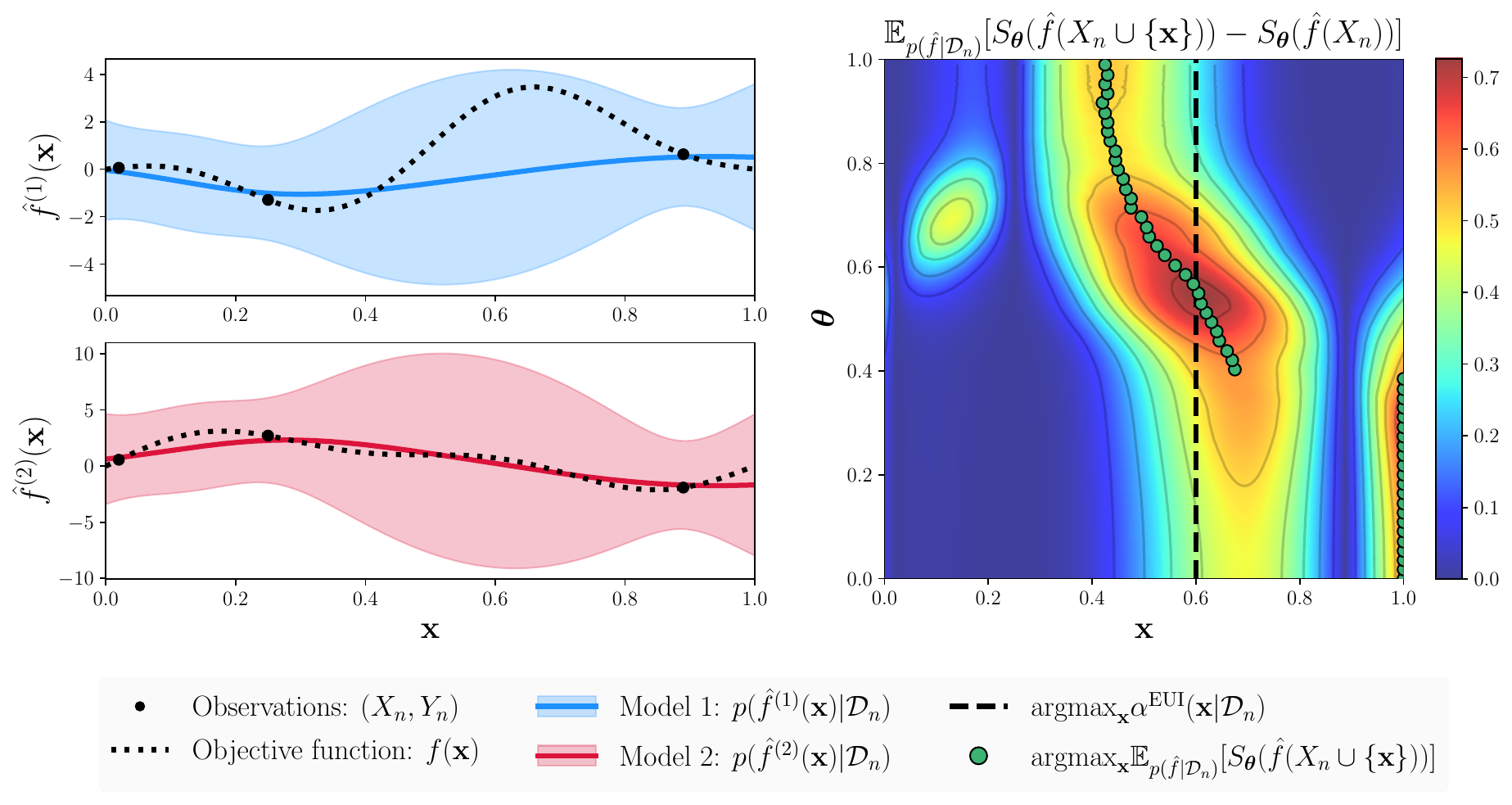}
	\centering
	\caption{A visualisation of the EUI acquisition function on a toy problem.}
	\label{fig:eui_acq_function}
\end{figure}

\begin{figure}
	\includegraphics[width=1\linewidth]{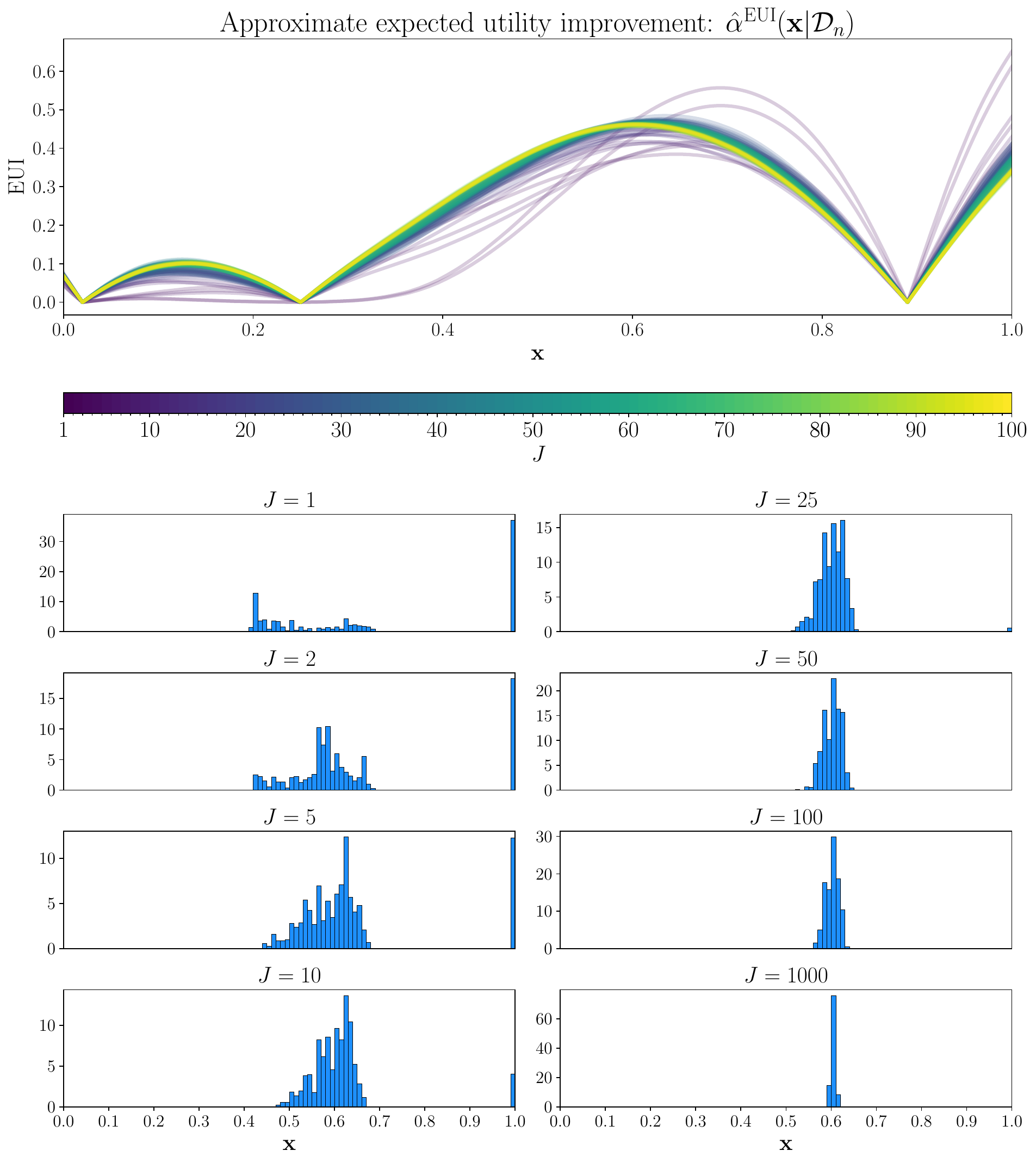}
	\centering
	\caption{A visualisation of the approximate EUI acquisition function with different numbers of Monte Carlo samples $J$ on the example described in \cref{fig:eui_acq_function}. On the top, we present estimates of the EUI acquisition function obtained using different numbers of Monte Carlo samples---the initial random seed is the same for each estimate. On the bottom, we present the density plots of the maximisers obtained from 1000 independent runs.}
	\label{fig:approximate_greedy_acq_density}
\end{figure}

\begin{example}
	[Toy problem] In \cref{fig:eui_acq_function}, we illustrate the approximate EUI acquisition function \eqref{eqn:approximate_eui} on a simple toy problem with two objectives defined on the unit interval $\mathbb{X} = [0,1]$. We consider the hypervolume indicator \eqref{eqn:hypervolume_indicator} as our R2 utility and assume a  Gaussian process prior on the function with a Gaussian observation likelihood. On the left, we present the mean and 99\% credible interval of the Gaussian process posteriors. On the right, we present the contours of the expected scalarised improvement over both the input space and scalarisation parameter space. The black dotted line highlights the best solution according to the EUI acquisition function, whilst the green dots highlights the maximisers for the expected scalarised improvement for each scalarisation parameter. 
	
	In \cref{fig:approximate_greedy_acq_density}, we visualise the approximate EUI acquisition function over many different samples. We see that the maximiser of the approximation slowly converges to the maximiser of the exact acquisition function when the number of samples increases. 
	\label{eg:toy_problem}
\end{example}
\subsubsection{Special cases}
\label{sec:discussion}

The family of AEUI acquisition functions \eqref{eqn:aeui} contains many existing multi-objective acquisition functions as special cases. In what follows, we will briefly review these different variants. Loosely speaking, the main differences among these approaches come from the choice of R2 utility, the number of Monte Carlo samples and the surrogate model.

\paragraph{Expected improvement} The EUI acquisition function \eqref{eqn:eui} can be interpreted as the multi-objective generalisation of the traditional expected improvement criterion \cite{mockus1975otitcn,jones1998jogo}. Many researchers have previously considered using this acquisition function before for some specific choices of utility functions \cite[Section 4]{zhan2020jgo}. For example, researchers have considered using the standard R2 utility \cite{deutz2019emo}, the hypervolume indicator \cite{emmerich2006itec}, and the IGD indicator \cite{astudillo2020icais}---all of which are special cases of the R2 utilities (\cref{sec:special_cases}). Following the example of Deutz, Emmerich, and Yang \cite{deutz2019emo} and Astudillo and Frazier \cite{astudillo2020icais}, we focussed our attention on the general setting, in which the utility is taken to be any arbitrary R2 utility. As illustrated in both of these papers, the EUI acquisition function\footnote{In these works, the EUI acquisition function with the R2 utility was referred to as the expected R2 improvement (ER2I) \cite{deutz2019emo} and the expected improvement under utility uncertainty (EI-UU) \cite{astudillo2020icais}, respectively.} can be efficiently estimated using Monte Carlo samples. Complementary to these works, we proved a general performance bound for this strategy in \cref{thm:aeui_guarantee}, which is to the best of our knowledge, is a new result. This finding is however unsurprising but it gives some justification for why these methods have been so successful in practice.

\paragraph{Thompson sampling} The multi-objective generalisation of the Thompson sampling \cite{thompson1933b} strategy can be recovered by using only a single model sample, $H=1$, in the approximate EUI acquisition function \eqref{eqn:approximate_eui}. This strategy has been proposed before for some specific utility functions. For instance, the TSEMO algorithm \cite{bradford2018jgo} considered the case where the hypervolume indicator is used as the utility.

\paragraph{Upper confidence bound} We obtain a multi-objective generalisation of the upper confidence bound (UCB) \cite{srinivas2010icml} strategy when we use a single function estimate, $\hat{f}: \mathbb{X} \rightarrow \mathbb{R}$, as our surrogate model. Ideally, this function estimate should be an upper bound for the true objective function. For instance, in standard Gaussian process regression, the function estimate could take the form $\hat{f}^{(m)}(\mathbf{x}) = \mu_n^{(m)}(\mathbf{x}) + \beta (\Sigma_n^{(m)}(\mathbf{x}, \mathbf{x}))^{1/2}$, where $\mu^{(m)}_n: \mathbb{X} \rightarrow \mathbb{R}$ is the mean function and $\Sigma^{(m)}_n: \mathbb{X} \times \mathbb{X} \rightarrow \mathbb{R}$ is the covariance function at time $n$ for objectives $m=1,\dots,M$, whilst $\beta \geq 0$ is a trade-off parameter that controls the width of the bound. One example of this variant being used in practice is the SMS-EGO algorithm \cite{ponweiser2008ppsn-px}, which considers the case where the hypervolume indicator is used as the utility function.

\paragraph{Random scalarisation} We recover the random scalarisation strategy when we use only a single sample, $J=1$, to the estimate the utility in the approximate EUI acquisition function \eqref{eqn:approximate_eui}. One of the earliest examples of this approach being used is the popular ParEGO algorithm \cite{knowles2006itec}, which considered the specific setting where the utility was the augmented variant of the standard R2 utility \eqref{eqn:standard_r2_utility}. Since then, many other researchers have proposed different variations of this strategy. For example, Paria, Kandasamy, and Póczos \cite{paria2020uai} proposed the Thompson sampling and UCB variants of this strategy when the surrogate model is a Gaussian process. They showed that both of these approaches satisfy some Bayesian regret bounds when we additionally assume that the family of scalarisations functions are both monotonic and Lipschitz. Subsequently, Zhang and Golovin \cite{zhang2020icml}, specialised these results for the hypervolume indicator. Complementary to these works, Chowdhury and Gopalan \cite{chowdhury2021icais} proved some frequentist regret bounds for the UCB variant of the random scalarisation strategy based on the approximate AEUI acquisition function \eqref{eqn:approximate_aeui}.

\paragraph{Other acquisition functions} There exists several other types of acquisition functions for multi-objective Bayesian optimisations, such as information-theoretic acquisition functions \cite{garrido-merchan2019n,belakaria2019anips,tu2022anips}, game-theoretic acquisition functions \cite{picheny2019jgo, binois2020jmlr}, uncertainty-based acquisition functions \cite{picheny2015sc, zuluaga2016jmlr}, and many others \cite{konakoviclukovic2020anips, malkomes2021icml}. These acquisition functions are not necessarily an instances of the AEUI acquisition function, but they have been shown to be an effective strategy to solve the R2 utility optimisation problem \eqref{eqn:r2_utility_problem}, in some cases, as illustrated in the aforementioned references. Having said that, in this paper we present a general approach in which the same utility function is used for both the optimisation procedure and the performance assessment because this leads to more interpretable results and clearer conclusions. Nevertheless, there are merits to using the other acquisition functions instead.
\section{Numerical experiments}
\label{sec:experiments}
We will now present some case studies to assess the effectiveness of the Bayesian optimisation strategies described in \cref{sec:bayesian_optimisation} and the corresponding performance bound in \cref{thm:aeui_guarantee}. For convenience, we will consider a standard Bayesian optimisation set-up, where our surrogate model is a Gaussian process which is obtained by placing a Gaussian process prior on the objective function and a Gaussian likelihood on the observations. We have implemented all of the algorithms using the Python library BoTorch \cite{balandat2020anips}. For all of the experiments, we initialised the data set by evaluating the function on $2(D+1)$ inputs which were obtained from uniformly sampling over the input space. For all of the performance plots, we have displayed the mean and two standard deviation of the log utility regret obtained from $100$ independent runs. The code that was used to run the experiments is available in our Github repository: \href{https://github.com/benmltu/scalarize}{\texttt{https://github.com/benmltu/scalarize}}. 

\begin{figure}
	\includegraphics[width=1\linewidth]{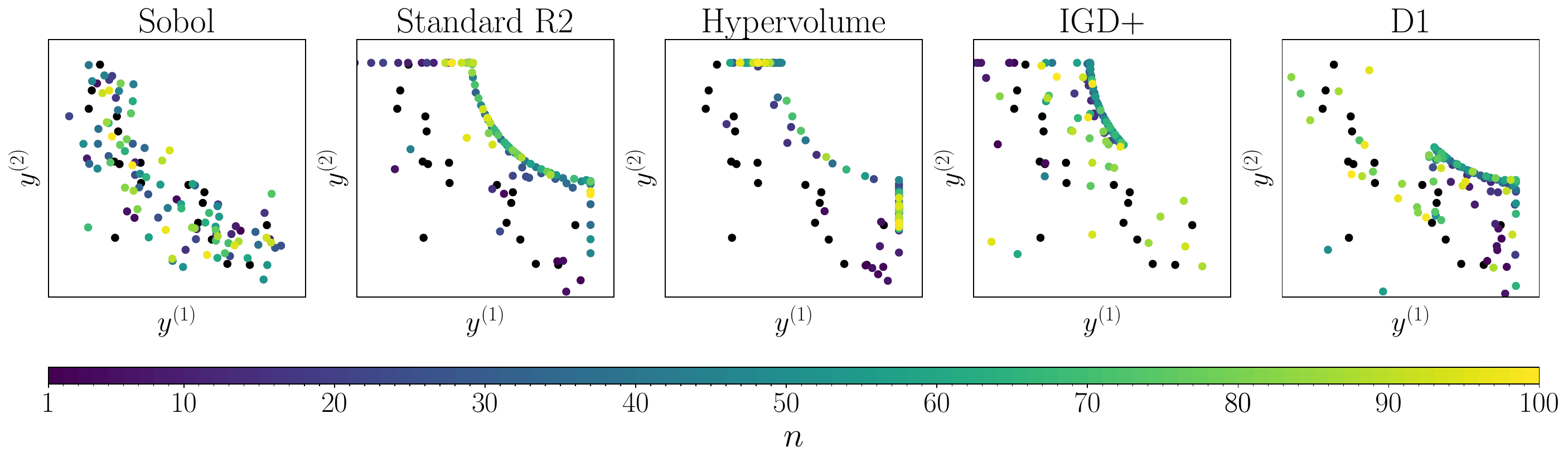}
	\centering
	\caption{A visualisation of the vectors obtained after one run of the approximate EUI approach with different choices of R2 utilities on the DTLZ2 benchmark with $D=8$ and $M=2$. We set $J=256$ and $H=128$ for the EUI estimate in \eqref{eqn:approximate_eui}.}
	\label{fig:dtlz2_output_space}
\end{figure}
\subsection{Different utilities}
\label{sec:different_utilities}
The performance of the EUI acquisition function is heavily dependent on the choice of performance metric. To illustrate this, we consider optimising the bi-objective DTLZ2 benchmark \cite{deb2002p2cecccn} using the four different R2 utilities outlined in \cref{sec:special_cases}. Specifically, we use the standard R2 utility and hypervolume indicator  in order to assess the whole Pareto front, but the IGD+ and D1 utilities are used to assess the upper and lower part of the Pareto front, respectively. We compare all of these greedy approaches with the random search baseline: Sobol.

In \cref{fig:dtlz2_output_space}, we present a visualisation of the objective vectors obtained after one run of Bayesian optimisation using the approximate EUI acquisition functions on the DTLZ2 benchmark. We see clearly that the choice of utility function influences the points that are selected. If we use a performance metric which assesses the whole Pareto front, then we end up querying points that try to cover the whole Pareto front. In contrast, if we use a performance metric which targets a particular region of the Pareto front, then we end up with more evaluations in this region.

\begin{figure}
	\includegraphics[width=1\linewidth]{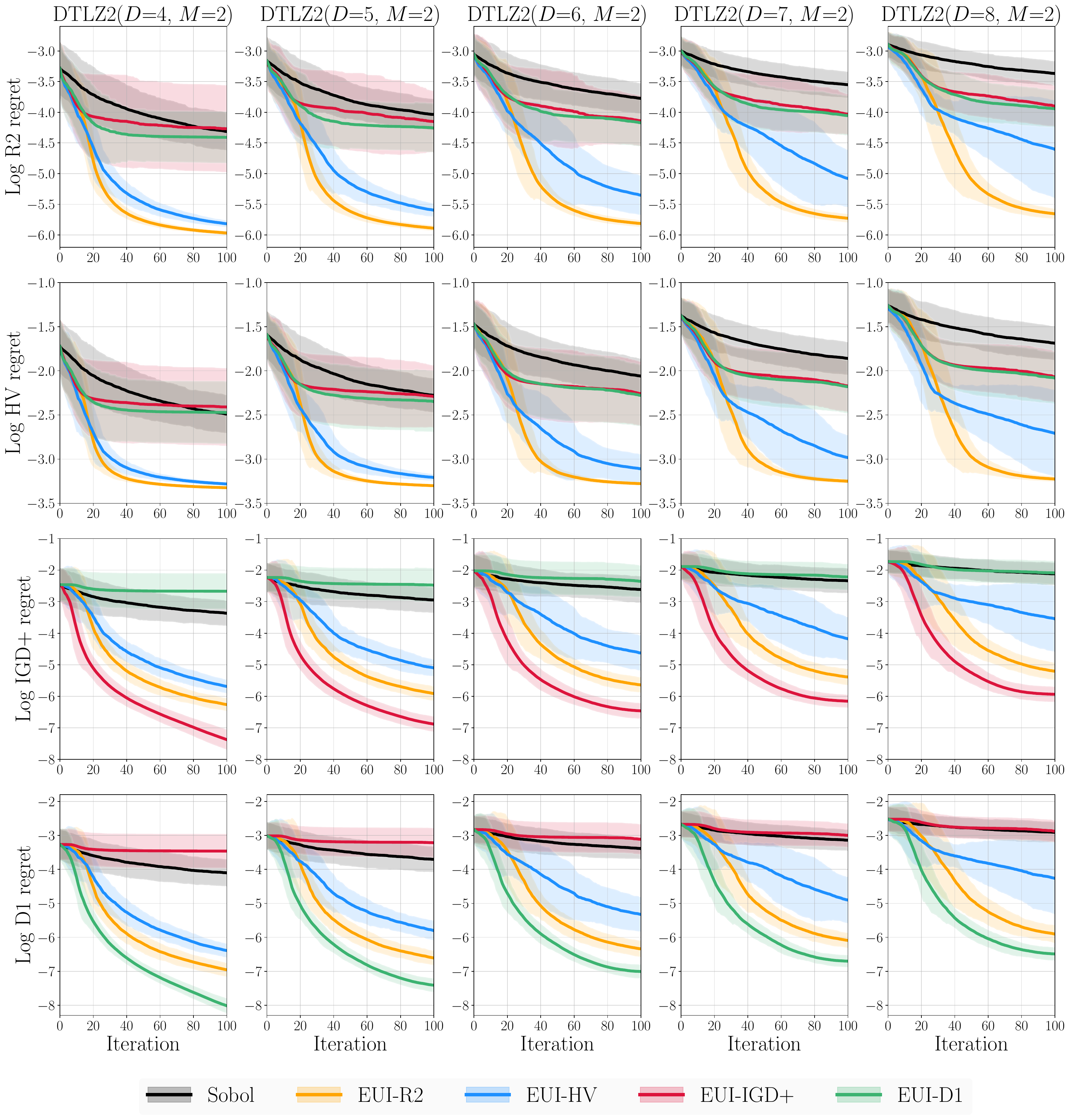}
	\centering
	\caption{A performance comparison of the approximate EUI approach with different choices of R2 utilities when applied on the DTLZ2 problem. In each row, we consider a different R2 utility, whereas in each column we consider a different number of inputs. For the EUI estimate \eqref{eqn:approximate_eui}, we set $J=256$ and $H=128$.}
	\label{fig:dtlz2_performance}
\end{figure}

In \cref{fig:dtlz2_performance}, we plot the performance of these acquisition functions over the different performance metrics. On the whole, we see that the performance of the EUI acquisition function is either the best or close to it on the performance metric that it was designed on. This pattern continues even after we increase the dimensionality of the input space. Interestingly, there are some cases in which we achieve better performance on one utility by using a different utility in the acquisition function. For example, we noticed that the EUI-R2 approach performed slightly better than the EUI-HV in terms of the hypervolume. Visually, the EUI-R2 approach manages to identify many more central points than the EUI-HV approach. This result is likely a consequence of the geometry of the Pareto front and the choice of reference points. In particular, after normalising the objectives to lie in the unit hypercube, the standard R2 utility was defined as an average weighted distance to the reference point $\boldsymbol{\upsilon}=(1.1, 1.1)$, whereas the hypervolume indicator was defined as an average weighted distance away from the reference point $\boldsymbol{\eta}=(-0.1, -0.1)$. As the Pareto front in this problem is concave, this resulted in superior performance for the former approach because it placed more emphasis on the central points. In contrast, the latter approach placed more emphasis on the corners and therefore had less coverage. 

For some additional intuition on the effect of changing the utility, we have included some additional examples in \cref{app:additional_figures} which studies the case where the Pareto front is three-dimensional.
\subsection{Different configurations}
\label{sec:different_configurations}
The performance bound in \eqref{eqn:aeui_guarantee} indicates that the quality of the approximate EUI acquisition function is dependent on the number of Monte Carlo samples and the quality of the model. To test this out, we analysed the empirical performance of this method on four different benchmark problems, which have input dimensions $D\in\{4,5,6,7\}$ and output dimensions $M\in\{3,4\}$. We used the implementation of these problems provided in the following references \cite{tanabe2020asc,daulton2021anips}.
\begin{itemize}
	\item \textbf{Rocket injector.} The rocket injector problem \cite{vaidyanathan20034asmae} considers optimising the performance and efficiency of a component used in a rocket. The objective function is defined over a $(D=4)$-dimensional space and consists of $M=3$ objectives.
	\item \textbf{Vehicle safety.} The vehicle design problem \cite{youn2004smo} considers optimising some variables relating to the safety of a vehicle. This objective function is defined over a $(D=5)$-dimensional space and consists of $M=3$ objectives.
	\item \textbf{Bulk carrier.} The bulk carrier design problem \cite{parsons2004josr} considers optimising some variables relating to the design specifications of a bulk carrier. This objective function is defined over a $(D=6)$-dimensional space and consists of $M=4$ objectives. 
	\item \textbf{Car side impact.} The car side-impact problem \cite{deb2009itec} considers optimising the performance of a car subject to some safety regulations. This objective function is defined over a $(D=7)$-dimensional space and consists of $M=4$ objectives.
\end{itemize}
\begin{figure}
	\includegraphics[width=0.9\linewidth]{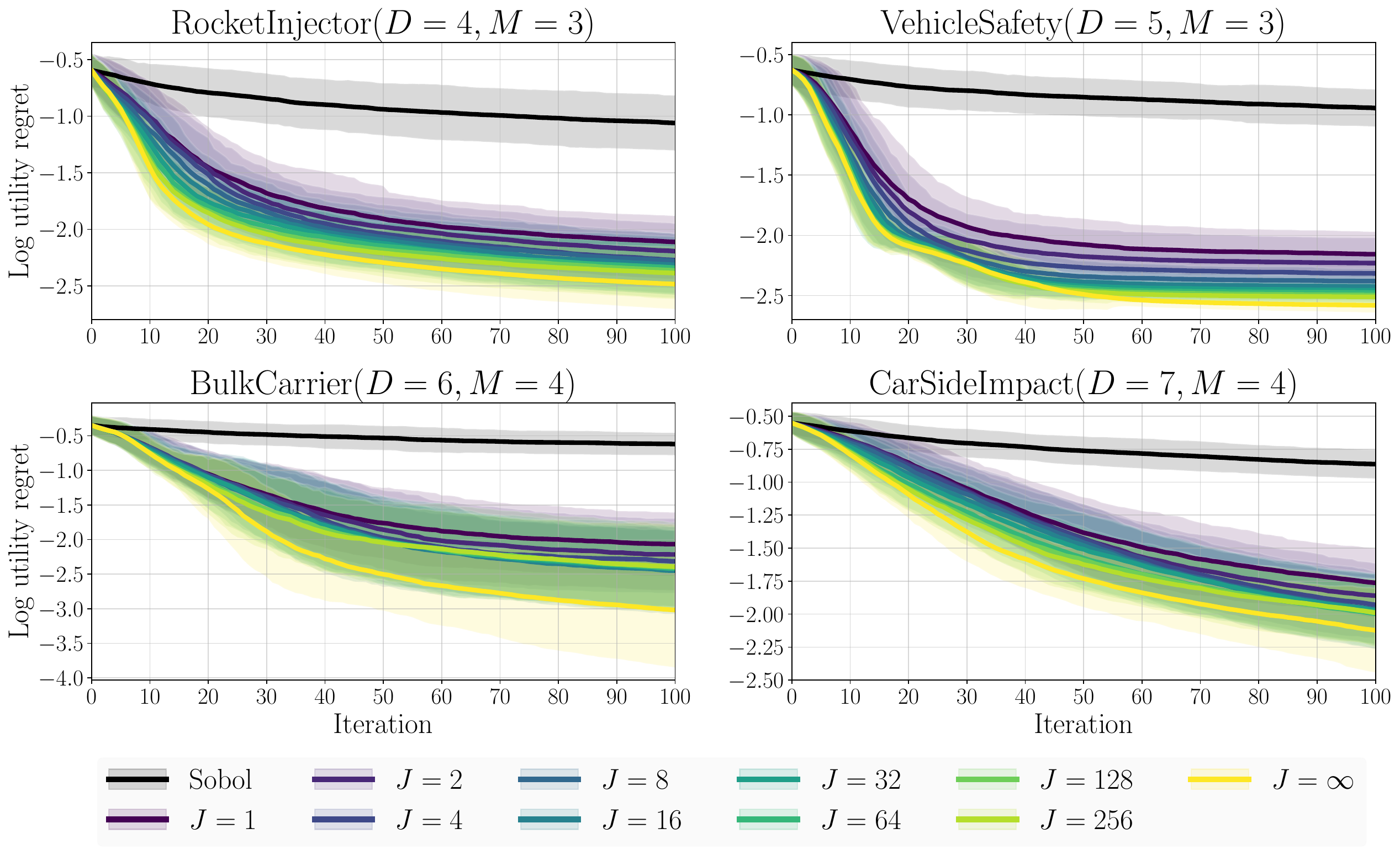}
	\centering
	\caption{A performance comparison of the approximate EUI approach when varying the number of samples used to estimate the utility. The utility function is set to the hypervolume indicator, whilst the number of model samples is set to $H=128$.}
	\label{fig:mc_performance_utility}
\end{figure}

In \cref{fig:mc_performance_utility}, we plot the performance of the approximate EUI acquisition function when we vary the number of scalarisation parameter samples. The label ``$J=\infty$'', denotes the setting where we compute the hypervolume indicator exactly \cite{daulton2021anips}. On the whole, we see that the performance does indeed improve when we increase the number of samples. This improvement is, however, quite small, which indicates that the variance of the Monte Carlo error is rather minor for these problems.
\begin{figure}
	\includegraphics[width=0.9\linewidth]{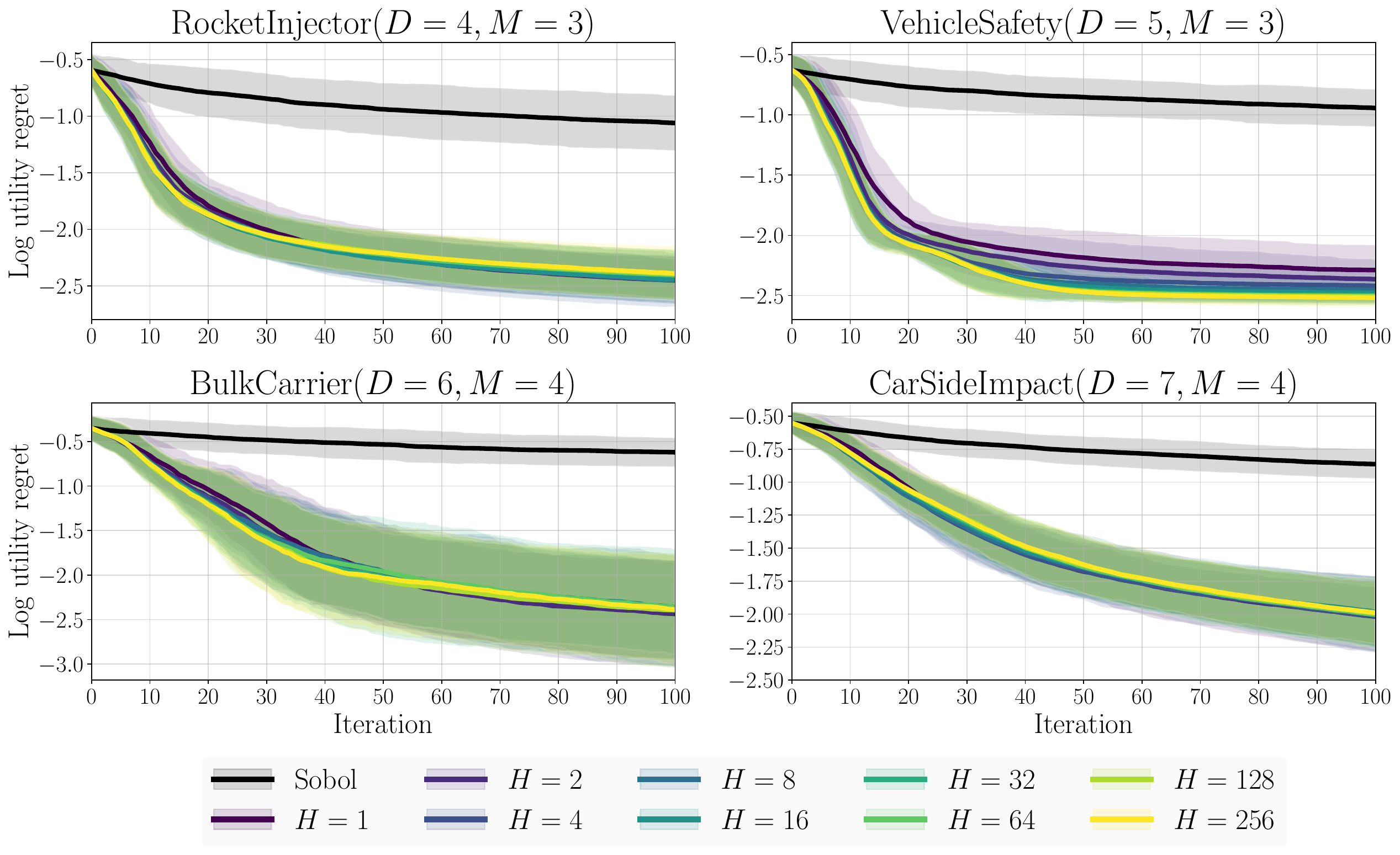}
	\centering
	\caption{A performance comparison of the approximate EUI approach when varying number of samples used to estimate the expectation over the model. The utility function is set to the hypervolume indicator, whilst the number of scalarisation parameter samples is set to $J=256$.}
	\label{fig:mc_performance_model}
\end{figure}

In \cref{fig:mc_performance_model}, we plot the performance of the approximate EUI acquisition function when we vary the number of model samples. We see that the difference in performance is very minor in most problems. There does not appear to be much computational benefit from increasing the number of model samples. In many cases, the Thompson sampling approach, where we just use one sample, performs just as well as the more expensive approaches. 
\begin{figure}
	\includegraphics[width=0.9\linewidth]{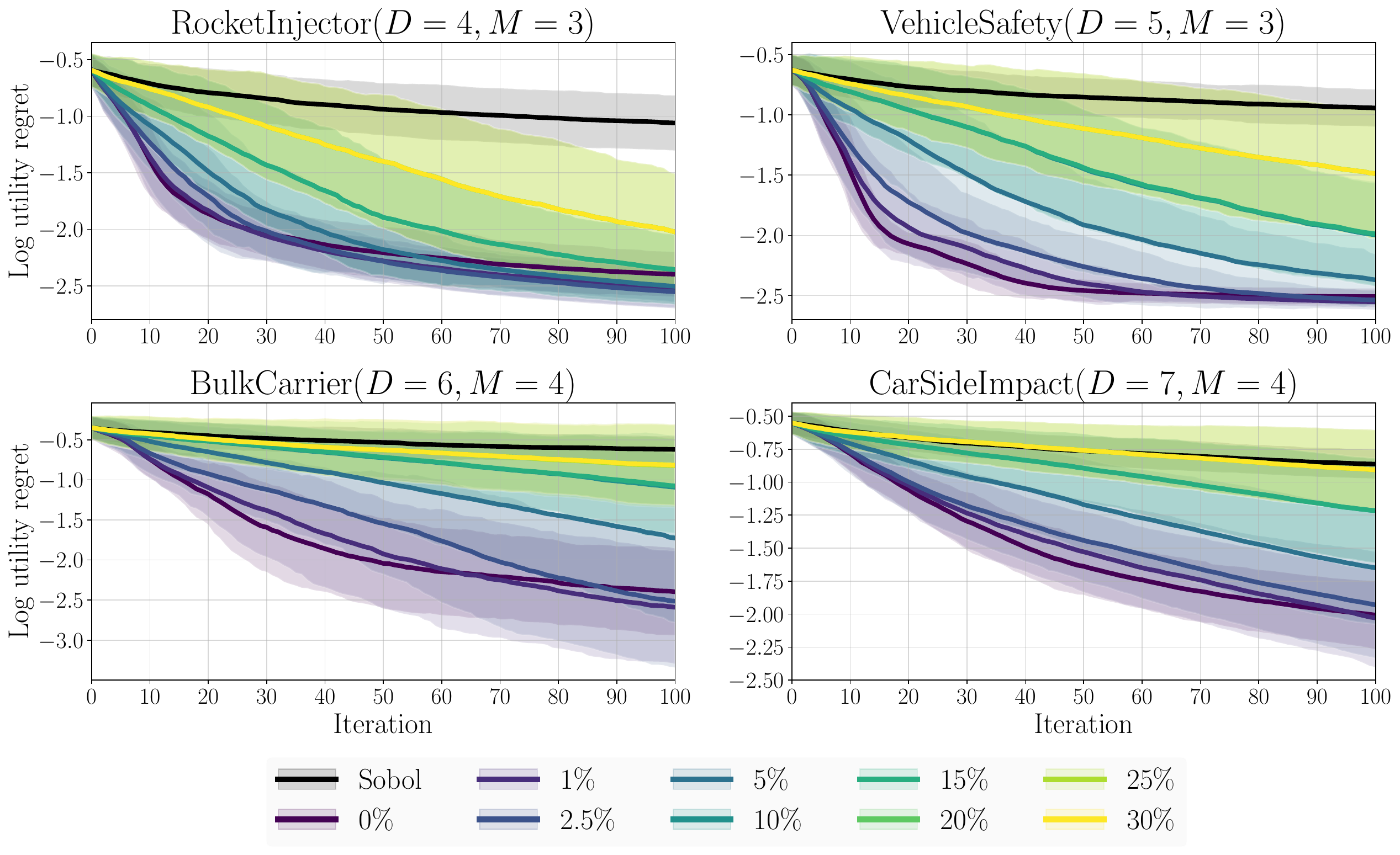}
	\centering
	\caption{A performance comparison of the approximate EUI approach when varying the standard deviation of the additive Gaussian noise. The standard deviation is set at a percentage of the objective ranges. The utility function is set to hypervolume indicator, whilst the number of Monte Carlo samples is set to $J=256$ and $H=128$.}
	\label{fig:noise_performance}
\end{figure}

In \cref{fig:noise_performance}, we plot the performance of the approximate EUI acquisition function when we vary the standard deviation of the additive Gaussian noise which we add to the function evaluations. As expected, the performance degrades when we incorporate more noise into the data. Surprisingly, there are some instances in which adding a small bit of noise actually improves the performance very slightly over the zero noise setting. This is perhaps a consequence of some additional exploration that takes place when the observations are a bit noisy.
\subsection{Mixed approaches}
\label{sec:mixed_approaches}
\begin{figure}
	\includegraphics[width=1\linewidth]{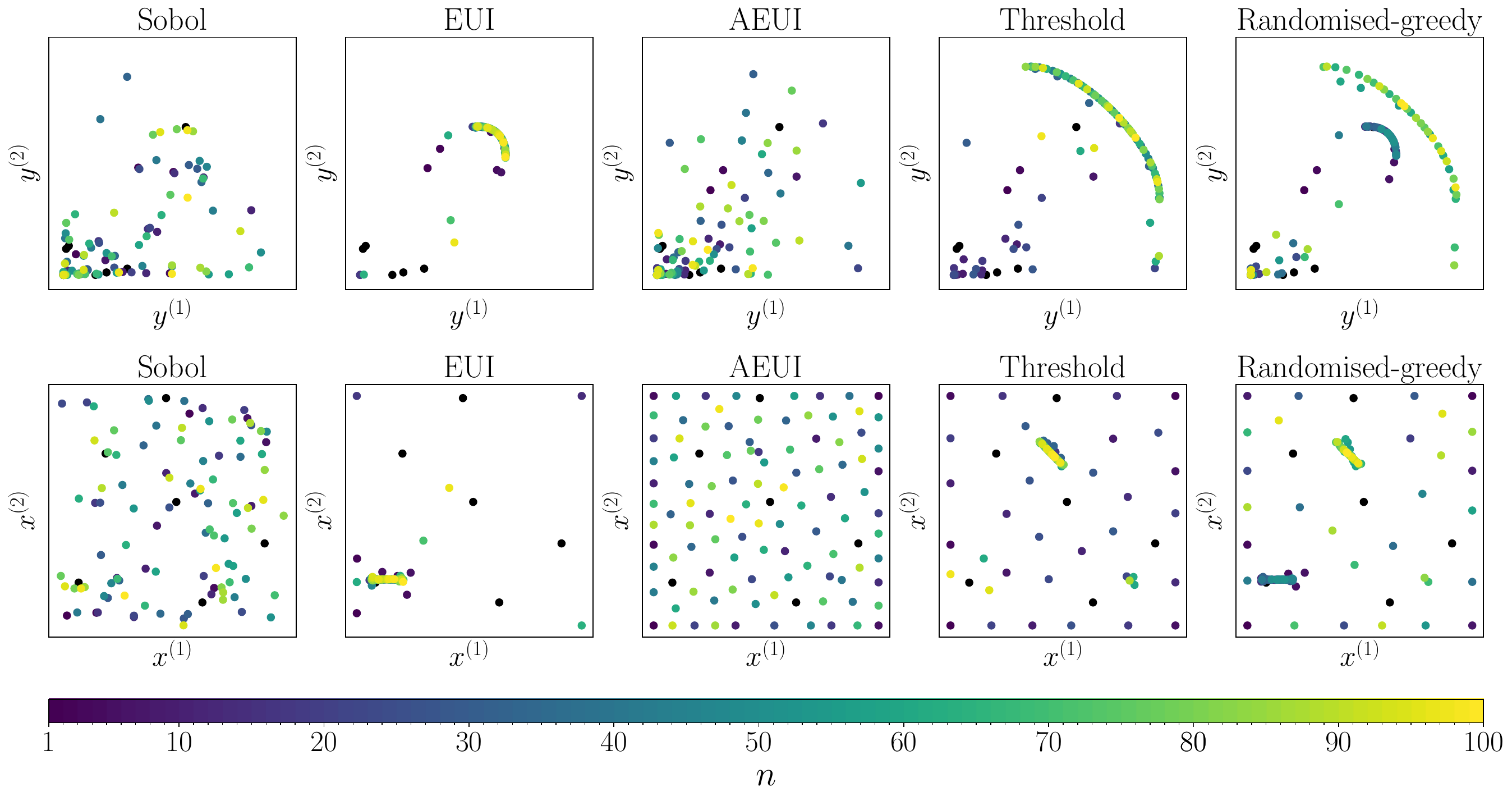}
	\centering
	\caption{A visualisation of the points obtained after one run of the different algorithms applied on the bi-objective GMM problem with $D=2$. For the mixed approaches, we set $p=0.3$.}
	\label{fig:gmm2_output_space}
\end{figure}

The EUI acquisition function can work well when the surrogate model is a good representation of the underlying objective. On other hand, there are some cases in which the EUI acquisition function can be overly greedy in its evaluations and can end up getting stuck in a suboptimal region. To illustrate this, consider the Gaussian mixture model (GMM) problem \cite{daulton2022icml}, where each objective is a mixture of Gaussian probability densities. This problem is tricky because the objective function has multiple modes which could potentially mislead the EUI acquisition function. We can see an instance of this happening in the second column of \cref{fig:gmm2_output_space}; specifically, we see that the EUI strategy manages to quickly locate one of the suboptimal modes and then only queries points near this region for the remaining iterations.

One strategy to overcome the over-exploitative behaviour of the EUI approach is to include a non-constant adjustment term into the acquisition function in order to encourage additional exploration. In the third column of \cref{fig:gmm2_output_space}, we illustrate the performance of this AEUI strategy when we use the adjustment term in \eqref{eqn:gp_adjustment} with $b_\delta = 1/M$. Unlike the EUI approach, we do not get stuck in one of the modes. On the contrary, we end up mimicking a space-filling design, where most of the function evaluations are spent in exploration in order to improve the model; only a few evaluations are actually spent in exploiting the model. In order to strike a better balance between exploration and exploitation, we propose using a novel mixed approach where we alternate between an explorative evaluation and an exploitative one. More precisely, we consider the mixed acquisition function
\begin{equation*}
	\alpha^{\text{EUI}-\lambda_n}(\mathbf{x}|\mathcal{D}_n) 
	= \alpha^{\text{EUI}}(\mathbf{x}|\mathcal{D}_n) + \lambda_n \mathcal{A}_{\delta}(\mathbf{x}|\mathcal{D}_n),
\end{equation*}
where $\lambda_n \geq 0$ is a time-dependent variable. For example, we could use the \emph{threshold approach}, where we set the schedule deterministically: $\lambda_n = \mathbbm{1}[n \leq p N]$ for some $p\in[0,1]$. Alternatively, we could use a \emph{randomised-greedy}\footnote{The randomised-greedy approach bares some similarity to the $\epsilon$-greedy approach which was previously proposed for the single-objective expected improvement \cite{death2021atelo}. In the $\epsilon$-greedy approach, we sample uniformly from the input space with probability $\epsilon \in [0, 1]$ as opposed to optimising a different acquisition function.} strategy, where we set the schedule randomly by sampling from a Bernoulli distribution: $\lambda_n \sim \text{Bernoulli}(p)$ for some $p\in[0,1]$.  As we can see in the final two columns of \cref{fig:gmm2_output_space}, the mixed approach manages to find the best of both worlds by spending some evaluations in improving the model and the other evaluations in exploiting what has been learnt.
\begin{figure}
	\includegraphics[width=1\linewidth]{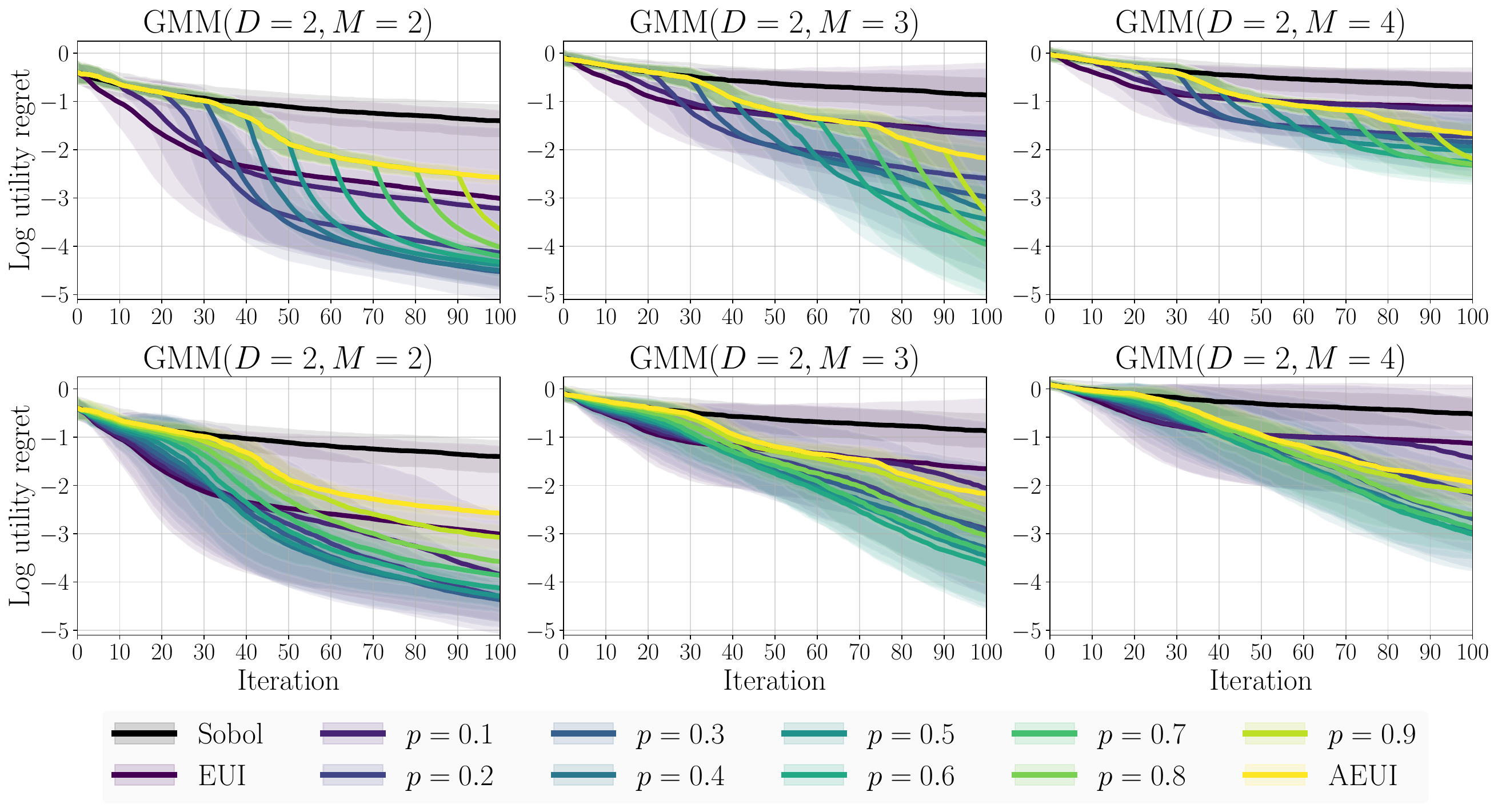}
	\centering
	\caption{A performance comparison of the threshold approach (top row) and randomised-greedy approach (bottom row) for different values of $p \in [0, 1]$ on the GMM problem with different number of objectives. The utility function was set to the hypervolume indicator.}
	\label{fig:gmm_performance}
\end{figure}
\begin{remark}
	[Mixed performance guarantee] Note that the performance bound for the mixed approach is justa  convex combination of the original bounds in \eqref{eqn:aeui_guarantee}. In particular, the Monte Carlo error remains the same, whereas the model error is a convex combination of the EUI model error and the AEUI model error; that is, $\epsilon_2 = (1-p) \epsilon^\textnormal{EUI}_2 + p \epsilon^\textnormal{AEUI}_2$.
\end{remark}

In \cref{fig:gmm_performance}, we present the performance plots of the mixed approaches on the GMM problem when we vary the number of objectives. We plot the logarithm of the hypervolume regret over many possible choices of $p \in [0, 1]$. Overall, we see that the performance of the mixed approach improves with the parameter $p$ up to a specific point before it begins to deteriorate again. This result perhaps gives some indication that there is an `optimal' choice of $p$, which likely depends on the utility, the surrogate model and the unknown objective function. In practice, we believe that this parameter should be set by the decision maker based on their confidence in the surrogate model. Specifically, a practitioner can use this parameter to modulate the number of evaluations which should be spent in exploration in order to improve the model and the number that should be spent in exploitation in order to improve the utility. 
\section{Summary and future work}
\label{sec:future_work}
In this paper, we presented a unified view of multi-objective optimisation. We described the two main philosophies that are used to solve this problem and then showed how they can be unified via the R2 utility functions. We then discussed the practical problem of actually solving this set optimisation problem \eqref{eqn:r2_utility_problem} once we have identified a suitable R2 utility. We focussed our discussion on the use of greedy algorithms because they come with theoretical performance guarantees. Among other things, we analysed the effectiveness of these greedy strategies in the popular Bayesian optimisation setting, where the objective function is modelled using a surrogate model. Notably, the analysis we have provided supports many of the Bayesian optimisation strategies that have been proposed over the past decade. To support these theoretical results, we also conducted several numerical experiments which demonstrate the effectiveness of our general methodology in practice. During our exposition, we also showcased some new variations of the greedy algorithms which arise naturally from our general framework. Overall, this paper has tried to consolidate many different ideas in multi-objective optimisation into one common framework. In what follows we discuss how some of these ideas could be used as a basis for future work.
\subsection{Preference elicitation}
\label{sec:preference_elicitation}
Practitioners often resort to the use of off-the-shelf performance metrics in order to assess the quality of a Pareto front approximation. As described in \cref{sec:r2_discussion}, this is counter-productive in the setting where a decision maker has preferences which are not adequately encoded into these performance metrics. The R2 utilities offers a potential solution to this problem. By carefully choosing the family of scalarisation function and parameter distributions, we can describe very complex preferences in the objective space. A useful direction for future work would focus on giving some practical guidance for setting the R2 utility in many commonly occurring situations. Alternatively, one might consider designing new adaptive strategies which try to learn the utility directly from user-interactions---this line of inquiry is related to the topic of preference learning. For example, the papers by Zintgraf et al. \cite{zintgraf2018p1icaams} and Lin et al. \cite{lin2022icais} considered the special case where the decision maker's preferences are described completely by a single scalarisation function. They then proceeded to learn this scalarisation function using different types of preference information. Notably, both of these papers were largely inspired by the Bayesian preference learning strategy of Chu and Ghahramani \cite{chu2005icml}. In this approach, a Gaussian process prior is placed on the scalarisation function and a likelihood is placed on the preference data. The scalarisation function is then updated using a Bayesian posterior update whenever more preference information is presented. For instance, Lin et al. \cite{lin2022icais} studied the setting in which the preference data came in the form of pairwise comparisons which were modelled using a probit likelihood, whereas Zintgraf et al. \cite{zintgraf2018p1icaams} also considered other preferences as well such as ranking and clustering preferences, which were modelled accordingly. 
\subsection{Constrained optimisation}
\label{sec:constrained_optimisation}
Real-world problems typically contain constraints which must also be taken into account. In this paper, we have largely overlooked this element of the problem and have implicitly absorbed any notion of a constraint into the definition of the input space $\mathbb{X} \subseteq \mathbb{R}^D$. Naturally, many of the general ideas and results that we have discussed here remains the same when we include constraints into the different scalar-valued optimisation problems such as in \eqref{eqn:scalarised_problem} and \eqref{eqn:utility_problem}. Nevertheless, there are also many practical considerations that have to be considered when we actually want to solve these constrained problems in practice. For instance, if we wanted to employ one of the greedy strategies described in \cref{sec:optimisation_problem}, then we would need to be able to optimise over this constrained space. This could perhaps be accomplished using a penalty-based optimisation strategy, although there are clearly many other possibilities as well \cite{nocedal2006}. These extensions are in general non-trivial and requires some additional thought. For example, with a penalisation strategy one has to decide what penalisation to use and where to incorporate it into the problem. Should we include it at the level of the objective function or at the level of the scalarisation or utility function? All of these choices lead to different algorithms and issues which would be interesting to study in their own right.

\subsection{Robustness}
\label{sec:robust_optimisation}
In real-world problems, the objective function is subject to many sources of uncertainty. A natural framework in which to study these types of problems is robust optimisation \cite{ben-tal2009}. Some effort has been spent over the last few decades in developing these robustness concepts for the multi-objective setting. Many of these extensions rely on the use of scalarisation functions \cite{ehrgott2014ejoor,ide2014mmor,fliege2014ejoor,daulton2022icml}. An interesting concept for future work would be to somehow unify many of these ideas into a common framework such as the one we presented in \cref{sec:r2_utility} for the standard setting.  
\subsection{Partially ordered optimisation}
The work in this paper has focussed on the multi-objective optimisation problem defined using the Pareto partial ordering. Nonetheless, many of the results that we presented here can also be extended to a more general setting in which the objective function $f: \mathbb{X} \rightarrow \mathbb{S}$, takes values in any partially ordered space $\mathbb{S}$. Specifically, we can follow the same recipe and reformulate the optimisation problem by appealing to the use of scalarisation functions $s: \mathbb{S} \rightarrow \mathbb{R}$, or utility functions $U: \mathbb{B}(\mathbb{S}) \rightarrow \mathbb{R}$. The R2 utilities can also be formulated and optimised in the same way as before. An interesting direction for future work could be focussed on investigating specific instances of this generalisation. For example, some researchers have previously attempted to incorporate user preferences (\cref{sec:preference_elicitation}) by directly defining some preference-based partial orderings \cite{drechsler1999ci,parmee2000ec}. Similarly, we can also treat the constrained (\cref{sec:constrained_optimisation}) and robust (\cref{sec:robust_optimisation}) optimisation problem as instances of this generalisation, where we appeal to the constrained \cite{fonseca1998itsmc-psh,deb2002itec} and robust \cite{ide2014mmor,ide2016os} Pareto domination rules instead.
\subsection{Decoupled evaluations}
We assumed throughout this paper that we always make full evaluations of the objective function. In some settings, we can choose to evaluate the objective function on only a subset of the objectives \cite{gelbart2014uai}. This might be advantageous when certain objectives are more expensive to evaluate than others. On the whole, it is not immediately clear how we can extend the formulation of the R2 utilities in order to handle this \emph{decoupled} set-up. The main difficulty lies in the fact that we have to assess the quality of a collection of partially observed objectives, which are not immediately comparable under the Pareto partial ordering. This decoupled framework is also related to optimisation problems for missing data \cite{bertsimas2018jmlr}. An interesting line of inquiry would be to focus on formalising this problem more concretely and then developing suitable algorithms to solve it.
\subsection{Multi-objective decision making}
Many of the ideas that we described here can also be exploited in other areas of research concerned with multi-objective decision making problems. Following Trivedi et al. \cite{trivedi2017itec}, we can broadly categorise multi-objective algorithms into three categories.
\begin{enumerate}
	\item \textbf{Domination-based.} These algorithms adopt the standard multi-objective optimisation perspective, where we solve the problem using only the sorting routines based on the Pareto partial ordering.
	\item\textbf{Decomposition-based.} These algorithms adopt the scalarisation perspective of multi-objective optimisation (\cref{sec:scalarisation_perspective}), which considers jointly solving a collection of single-objective problems using scalarisation functions.
	\item \textbf{Indicator-based.} These algorithms adopt the utility perspective of multi-objective optimisation (\cref{sec:utility_perspective}), which rely on the use of utility functions in order to determine which sets of vectors are the most desirable.
\end{enumerate}

In \cref{sec:r2_utility}, we demonstrated how all three categories are linked via the R2 utilities \eqref{eqn:r2_utility}. A consequence of this connection is that many indicator-based algorithms can potentially be converted into a decomposition-based algorithm via a suitably chosen distribution of scalarisation functions and vice versa. For example, we can exploit this connection in the development of new multi-objective evolutionary algorithms \cite{zhou2011saec, li2015acs, zhou2019}, multi-objective reinforcement learning algorithms \cite{roijers2013j, vanmoffaert20132isadprla, vanmoffaert2013emo, zintgraf2015amlcbn}, and multi-objective gradient-based optimisation algorithms (\cref{sec:gradient_based_optimisation}).
\subsection{Gradient-based optimisation}
\label{sec:gradient_based_optimisation}
The traditional goal of multi-objective gradient-based optimisation is to target a finite collection of Pareto stationary points \cite{kuhn1951potsbsomsap}. This is typically accomplished by appealing to the scalarisation methodology in which scalarised optimisation problems \eqref{eqn:scalarised_problem} are defined and consequently solved using standard gradient-based methods. There also exist many adaptive variants of these algorithms in which the scalarisation parameters do not have to be specified a priori, but can be updated in an online manner \cite{fliege2000mmoo,schaffler2002jootaa,fliege2009sjo,eichfelder2009sjo,desideri2012crm,sener2018anips,mahapatra2020icml}. Notably much less work has been focussed on solving the utility optimisation problem \eqref{eqn:utility_problem} with gradient-based methods \cite{emmerich2007hm,sosahernandez2014e-abbpsonaecv,wang2017emo,deist2021a}. The key challenge with this latter problem is that the (sub-)gradients\footnote{The gradient of a utility function is not necessarily well-defined in practice. For example, in an R2 utility \eqref{eqn:r2_utility} we would need to differentiate a maximum of some scalarised functions. To address this problem it is common to appeal to a more relaxed notion of a gradient \cite{parikh2014fto}.} of a Pareto compliant utility function are zero at any input which is Pareto dominated. In regards to an R2 utility \eqref{eqn:r2_utility}, this can happen when an input does not lead to a maximal scalarised value for any scalarisation parameter. This issue is clearly problematic in a gradient-based algorithm because whenever an input becomes sub-optimal, with respect to the current set, then it is no longer updated. Many authors have suggested including some additional dynamics to tackle this issue \cite{wang2017emo,deist2021a}. However, these methods have not been studied theoretically and their relationship to existing scalarisation-based algorithms is unclear. In particular, we anticipate that there are many useful ideas from the existing scalarisation-based algorithms that could be exploited in order to solve the utility optimisation problem in a more efficient and theoretically justifiable way. 

\section*{Acknowledgements}
Ben Tu was supported by the EPSRC StatML CDT programme EP/S023151/1 and BASF SE, Ludwigshafen am Rhein. Nikolas Kantas was partially funded by JPMorgan Chase \& Co. under J.P. Morgan A.I. Faculty Research Awards 2021. We would also like to acknowledge the associate editor and reviewers for their useful suggestions and insights that have helped to improve this paper.
\newpage
\appendix
\section{Proof of results}
\label{app:proofs}
\subsection{Proof of \cref{prop:monotonic_implies_optimal}} 
\label{app:proofs:prop:monotonic_implies_optimal}

These two results follow from a simple proof by contradiction. Let $s: \mathbb{R}^M \rightarrow \mathbb{R}$ denote a scalarisation function which is strictly monotonically increasing over the feasible objective space. Suppose that there exists a solution $\mathbf{x}^* \in X^*_s$ which is not strictly Pareto optimal. Then there must exist an input $\mathbf{x} \in \mathbb{X}$ such that $f(\mathbf{x}) \succ f(\mathbf{x}^*)$. By strict monotonicity, we have $s(f(\mathbf{x})) > s(f(\mathbf{x}^*))$, which contradicts the optimality of $\mathbf{x}^*$.

Similarly, for the other result, let $s: \mathbb{R}^M \rightarrow \mathbb{R}$ denote a scalarisation function which is strongly monotonically increasing over the feasible objective space. Suppose that there exists a solution $\mathbf{x}^* \in X^*_s$ which is not weakly Pareto optimal. Then there must exist an input $\mathbf{x} \in \mathbb{X}$ such that $f(\mathbf{x}) \succsucc f(\mathbf{x}^*)$. By strong monotonicity, we have $s(f(\mathbf{x})) > s(f(\mathbf{x}^*))$, which contradicts the optimality of $\mathbf{x}^*$.

\begin{flushright}
	$\blacksquare$
\end{flushright}

\subsection{Proof of \cref{prop:monotone_submodular}} 
\label{app:proofs:prop:monotone_submodular}
In the following paragraph, we will prove that the set function $S_{\boldsymbol{\theta}}$ satisfies the monotone and diminishing returns property for all scalarisation parameters $\boldsymbol{\theta} \in \Theta$. Once we have this result, we can conclude that any R2 utility, which is defined using these set functions, also satisfies these properties because inequalities are preserved when we compute the expectation with respect to a non-negative density $p(\boldsymbol{\theta}) \geq 0$.

Consider any finite set of vectors $A, B \in \mathbb{B}(\mathbb{R}^M)$ with $A \subseteq B$ and any vector $\mathbf{c} \in \mathbb{R}^M$. Then clearly the set function $S_{\boldsymbol{\theta}}$ is monotonic because $S_{\boldsymbol{\theta}}(A) \leq S_{\boldsymbol{\theta}}(B)$ for all $\boldsymbol{\theta} \in \Theta$. We will now show that $S_{\boldsymbol{\theta}}$ also satisfies the diminishing returns property. If the maximum values are the same, $S_{\boldsymbol{\theta}}(A) = S_{\boldsymbol{\theta}}(B)$, then the diminishing returns property is trivially satisfied. If the maximum values are different, $S_{\boldsymbol{\theta}}(A) < S_{\boldsymbol{\theta}}(B)$, then there are two settings to consider. In the first case, we have $S_{\boldsymbol{\theta}}(B) \geq S_{\boldsymbol{\theta}}(\{\mathbf{c}\})$, which implies $S_{\boldsymbol{\theta}}(A \cup \{\mathbf{c}\}) - S_{\boldsymbol{\theta}}(A) \geq 0 = S_{\boldsymbol{\theta}}(B \cup \{\mathbf{c}\}) - S_{\boldsymbol{\theta}}(B)$. Whilst in the second case, we have $S_{\boldsymbol{\theta}}(B) < S_{\boldsymbol{\theta}}(\{\mathbf{c}\})$, which implies
\begin{align*}
	S_{\boldsymbol{\theta}}(A \cup \{\mathbf{c}\}) - S_{\boldsymbol{\theta}}(A) 
	&= S_{\boldsymbol{\theta}}(\{\mathbf{c}\}) - S_{\boldsymbol{\theta}}(A) 
	\\
	&\geq S_{\boldsymbol{\theta}}(\{\mathbf{c}\}) - S_{\boldsymbol{\theta}}(B) 
	= S_{\boldsymbol{\theta}}(B \cup \{\mathbf{c}\}) - S_{\boldsymbol{\theta}}(B).
\end{align*}

\begin{flushright}
	$\blacksquare$
\end{flushright}

\subsection{Proof of \cref{prop:weakly_pareto_compliant}} 
\label{app:proofs:prop:weakly_pareto_compliant}
Consider any finite set of vectors $A, B \in \mathbb{B}(\mathbb{R}^M)$ with $A \succeq B$. As the scalarisation functions $s_{\boldsymbol{\theta}}$ are monotonically increasing for all scalarisation parameters, this implies $S_{\boldsymbol{\theta}}(A) \geq S_{\boldsymbol{\theta}}(B)$ for all $\boldsymbol{\theta} \in \Theta$. By computing the expectation with respect to the distribution $p(\boldsymbol{\theta}) \geq 0$, we obtain the final result $U(A) \geq U(B)$.

\begin{flushright}
	$\blacksquare$
\end{flushright}
\subsection{Proof of \cref{prop:hypervolume_indicator}} 
\label{app:proofs:prop:hypervolume_indicator}
The proof of this result is presented in the following references: \cite[Section 3.2]{shang2018pgecc}, \cite[Section 2]{deng2019itec} and \cite[Lemma 5]{zhang2020icml}.

\begin{flushright}
	$\blacksquare$
\end{flushright}
\subsection{Proof of \cref{prop:igd}} 
\label{app:proofs:prop:igd}
This result can be derived by simply expanding the definition of the IGD indicator \eqref{eqn:igd}:
\begin{align*}
	(I^{\text{IGD}_{p, q}}(Y, \Upsilon))^q 
	&= \frac{1}{|\Upsilon|} 
	\sum_{\boldsymbol{\upsilon} \in \Upsilon} \Bigl(
	\min_{\mathbf{y} \in Y} ||\boldsymbol{\upsilon} - \mathbf{y}||_{L^p}
	\Bigr)^q 
	\\
	&= \frac{1}{|\Upsilon|} 
	\sum_{\boldsymbol{\upsilon} \in \Upsilon}
	\min_{\mathbf{y} \in Y} ||\boldsymbol{\upsilon} - \mathbf{y}||^q_{L^p}
	\\
	&= - \frac{1}{|\Upsilon|} 
	\sum_{\boldsymbol{\upsilon} \in \Upsilon}
	\max_{\mathbf{y} \in Y} (-||\boldsymbol{\upsilon} - \mathbf{y}||^q_{L^p})
	\\
	&= - \frac{1}{|\Upsilon|} 
	\sum_{\boldsymbol{\upsilon} \in \Upsilon} \max_{\mathbf{y} \in Y} s^{\textnormal{IGD}_{p, q}}_{\boldsymbol{\upsilon}}(\mathbf{y})
\end{align*}
for any $Y, \Upsilon \in \mathbb{B}(\mathbb{R}^M)$ and $p, q \geq 1$.

\begin{flushright}
	$\blacksquare$
\end{flushright}
\subsection{Proof of \cref{prop:d1}} 
\label{app:proofs:prop:d1}
This result can be derived by simply expanding the definition of the D1 indicator \eqref{eqn:d1}:
\begin{align*}
	I^{\text{D1}}(Y, \Upsilon, \mathbf{w})
	&= \frac{1}{|\Upsilon|} \sum_{\boldsymbol{\upsilon} \in \Upsilon}
	\min_{\mathbf{y} \in Y} \max_{m=1,\dots,M} w^{(m)} (\upsilon^{(m)} - y^{(m)})
	\\
	&= - \frac{1}{|\Upsilon|} \sum_{\boldsymbol{\upsilon} \in \Upsilon}
	\max_{\mathbf{y} \in Y} \Bigl(- \max_{m=1,\dots,M} w^{(m)} (\upsilon^{(m)} - y^{(m)}) \Bigr)
	\\
	&= - \frac{1}{|\Upsilon|} \sum_{\boldsymbol{\upsilon} \in \Upsilon}
	\max_{\mathbf{y} \in Y} s^{\textnormal{Chb}}_{(\boldsymbol{\upsilon}, \mathbf{w})}(\mathbf{y})
\end{align*}
for any $Y, \Upsilon \in \mathbb{B}(\mathbb{R}^M)$ and $\mathbf{w} \in \Delta^{M-1}$.

\begin{flushright}
	$\blacksquare$
\end{flushright}
\subsection{Proof of \cref{thm:greedy_guarantee}} The proof of this result is identical to the one described by Krause and Golovin \cite[Theorem 1.5]{krause2014tpathp}. 

\begin{flushright}
	$\blacksquare$
\end{flushright}
\subsection{Performance bounds}
\label{app:proofs:performance_bounds}
The proof of our performance bounds follows a similar pattern to the proof of the original greedy guarantee described by Krause and Golovin \cite[Theorem 1.5]{krause2014tpathp}. To establish these probabilistic bounds, we will first prove a useful result (\cref{lemma:performance_bound}). For convenience of notation, in what follows we denote an optimal set of inputs by
\begin{equation*}
	X_P^* \in \argmax_{X \subseteq \mathbb{X}, |X|\leq P} U(f(X))
\end{equation*}
for any $P > 0$.
\begin{lemma}
	Consider an objective function $f: \mathbb{X} \rightarrow \mathbb{R}^M$, a non-negative R2 utility $U: \mathbb{B}(\mathbb{R}^M) \rightarrow \mathbb{R}_{\geq 0}$, and a collection of inputs and scalars $\{X_n\}_{n\geq1}$ and $\{\rho_n\}_{n\geq1}$, where $X_n \subseteq \mathbb{X}$ and $\rho_n \in \mathbb{R}$. Suppose for some positive integer $P$, that there exist a $\delta \in (0, 1)$ such that the following inequality holds with probability $1-\delta$ for $n \leq N$:
	\begin{equation*}
		U(f(X_P^*)) \leq U(f(X_n)) + P (U(f(X_{n+1})) - U(f(X_n))) + P \rho_{n+1}.
	\end{equation*}
	Then the following inequality holds with probability $1-N\delta$,
	\begin{equation*}
		U(f(X_N)) 
		\geq (1-e^{-N/P}) U(f(X_P^*)) - \sum_{n=1}^N \rho_n \Bigl(1-\frac{1}{P}\Bigr)^{N-n}.
	\end{equation*}
	\label[lemma]{lemma:performance_bound}
\end{lemma}
\paragraph{Proof of \cref{lemma:performance_bound}} Let $\zeta_n = U(f(X_P^*)) - U(f(X_n))$; then rearranging the inequality, we obtain
\begin{align*}
	\zeta_{n+1} \leq \Bigl(1-\frac{1}{P} \Bigr)\zeta_n + \rho_{n+1}.
\end{align*}
By repeated use of this inequality, we find that
\begin{align*}
	\zeta_N 
	&\leq \Bigl(1-\frac{1}{P}  \Bigr)^N \zeta_0 + \sum_{n=1}^N \rho_n \Bigl(1-\frac{1}{P} \Bigr)^{N-n},
\end{align*}
which is an inequality that holds with probability $1-N\delta$ by the union bound (Boole's inequality). As the R2 utility is non-negative, we have that $\zeta_0 = U(f(X_P^*))- U(f(\emptyset)) \leq U(f(X_P^*))$. Using this result and the fact that $1-x \leq e^{-x}$ for all $x \in \mathbb{R}$, we obtain the final performance bound
\begin{equation*}
	U(f(X_N)) 
	\geq (1-e^{-N/P}) U(f(X_P^*)) - \sum_{n=1}^N \rho_n \Bigl(1-\frac{1}{P} \Bigr)^{N-n}.
\end{equation*}

\begin{flushright}
	$\blacksquare$
\end{flushright}
\subsubsection{Proof of \cref{thm:approximate_greedy_guarantee}} 
\label{app:proofs:thm:approximate_greedy_guarantee}
We begin the proof by first establishing an inequality between the series of utility values:
\begin{align}
	U(f(X_P^*)) 
	&\leq U(f(X_P^* \cup X_n))
	\label{eqn:submodular_bound_1}
	\\
	&= U(f(X_n)) + \sum_{p=1}^P (U(f(X_p^* \cup X_n)) - U(f(X_{p-1}^* \cup X_n)))
	\label{eqn:submodular_bound_2}
	\\
	&\leq U(f(X_n)) + \sum_{\mathbf{x} \in X_P^*} (U(f(X_n \cup \{\mathbf{x}\})) - U(f(X_n)))
	\label{eqn:submodular_bound_3}
	\\
	&\leq U(f(X_n)) + P (U(f(X_n \cup \{\mathbf{x}^*\})) - U(f(X_n))).
	\label{eqn:submodular_bound_4}
\end{align}
The first line \eqref{eqn:submodular_bound_1} follows from the monotone property. The second line \eqref{eqn:submodular_bound_2} by from rewriting the utility as a telescopic sum. The third line \eqref{eqn:submodular_bound_3} by from applying the diminishing returns property and the last line \eqref{eqn:submodular_bound_4} by from setting the input $\mathbf{x} \in \mathbb{X}$ to be a maximiser, $\mathbf{x}^* \in \argmax_{\mathbf{x} \in X_P^*}(U(f(X_n \cup \{\mathbf{x}\})) - U(f(X_n)))$, and then using the fact that $|X_P^*| \leq P$, where $|\cdot|$ denotes the cardinality of the set.
\par
As the scalarisation parameters are sampled independently, we have the following result by Hoeffding's inequality:
\begin{align*}
	&\mathbb{P}\Bigl[\bigl|(\hat{U}_{J}(f(X_n \cup \{\mathbf{x}\})) - \hat{U}_{J}(f(X_n))) - 
	(U(f(X_n \cup \{\mathbf{x}\})) - U(f(X_n))) \bigr| 
	\geq \xi_n \Bigr] 
	\\
	&\leq 2 \exp\biggl( -\frac{2 J \xi_n^2}{C_n^2} \biggr)
\end{align*}
for any $\mathbf{x} \in \mathbb{X}$, where $C_n = \sup_{\boldsymbol{\theta} \in \Theta}(S_{\boldsymbol{\theta}}(f(\mathbb{X})) - S_{\boldsymbol{\theta}}(f(X_n)))$ for $n \geq 1$ and $C_0 := C$. By manipulating the expression, we find that this inequality will hold with probability $1-\delta$ if we set $\xi_n = C_n \sqrt{\log(2/\delta) / (2J)}$. Using this result, we obtain the following probability bound, which holds with probability $1-2\delta$ by the union bound:
\begin{align}
	U(f(X_P^*)) 
	&\leq U(f(X_n)) + P (
	\hat{U}_J(f(X_n \cup \{\mathbf{x}^*\})) - \hat{U}_J(f(X_n)) + \xi_n)
	\label{eqn:proof:approximate_greedy_a}	
	\\
	&\leq U(f(X_n)) + P(\hat{U}_J(f(X_{n+1})) - \hat{U}_J(f(X_n)) + \xi_n)
	\label{eqn:proof:approximate_greedy_b}	
	\\
	&\leq U(f(X_n)) + P (U(f(X_{n+1})) - U(f(X_n))) + 2P\xi_n.
	\label{eqn:proof:approximate_greedy_c}	
\end{align}
In the first equation \eqref{eqn:proof:approximate_greedy_a}, we used Hoeffding's inequality for $\mathbf{x}^* \in X^*_P$. In the second equation \eqref{eqn:proof:approximate_greedy_b}, we used the fact the input $\mathbf{x}_{n+1} \in \mathbb{X}$ is picked according to the approximate greedy strategy and in the last equation \eqref{eqn:proof:approximate_greedy_c}, we again used Hoeffding's inequality for the approximately greedily selected input $\mathbf{x}_{n+1} \in \mathbb{X}$. Applying \cref{lemma:performance_bound}, we find that the following inequality holds with probability $1-2N\delta$: 
\begin{equation*}
	U(f(X_N)) 
	\geq (1-e^{-N/P}) U(f(X_P^*)) - \sum_{n=1}^N 2 \xi_{n-1} \Bigl(1-\frac{1}{P} \Bigr)^{N-n}.
\end{equation*}

\begin{flushright}
	$\blacksquare$
\end{flushright}
\subsubsection{Proof of \cref{thm:aeui_guarantee}}
\label{app:proofs:thm:aeui_guarantee}
Firstly, by Hoeffding's inequality, we obtain the following concentration bound for the expectation over the surrogate model:
\begin{equation*}
	\mathbb{P}\biggl[
	\Bigl|\frac{1}{H}\sum_{h=1}^H \hat{U}_{J}(\hat{f}_h(X)) 
	-\mathbb{E}_{p(\hat{f}|\mathcal{D}_n)}[\hat{U}_{J}(\hat{f}(X))] \Bigr| \geq \xi_{n, 1} 
	\biggr] 
	\leq 2 \exp\biggl(
	-\frac{2 H \xi_{n, 1}^2}{C_n^2}
	\biggr)
\end{equation*}
where $X = X_n \cup \{\mathbf{x}\}$ for any $\mathbf{x} \in \mathbb{X}$ and $C_n = \sup_{\boldsymbol{\theta} \in \Theta}\mathbb{E}_{p(\hat{f}|\mathcal{D}_n)}[S_{\boldsymbol{\theta}}(\hat{f}(\mathbb{X})) - S_{\boldsymbol{\theta}}(\hat{f}(X_n))]$. By manipulating this expression, we find that the inequality holds with probability $1-\delta$ if we set $\xi_{n, 1} = C_n \sqrt{\log(2/\delta) / (2H)}$. Similarly, by applying Hoeffding's inequality again, we obtain the following concentration bound for the utility estimate:
\begin{equation*}
	\mathbb{P}\biggl[
	\Bigl|\mathbb{E}_{p(\hat{f}|\mathcal{D}_n)}[\hat{U}_{J}(\hat{f}(X))]
	-\mathbb{E}_{p(\hat{f}|\mathcal{D}_n)}[U(\hat{f}(X))] \Bigr| \geq \xi_{n, 2}
	\biggr] 
	\leq 2 \exp\biggl(
	-\frac{2 J \xi_{n, 2}^2}{C_n^2}
	\biggr).
\end{equation*}
By manipulating this expression, we find that the inequality holds with probability $1-\delta$ if we set $\xi_{n, 2} = C_n \sqrt{\log(2/\delta) / (2J)}$. We now upper bound the sum of these two parameters by
\begin{equation*}
	\xi_n := C_n \sqrt{\frac{2}{\min(J,H)} \log\biggl(\frac{2}{\delta}\biggr)} \geq \xi_{n, 1} + \xi_{n, 2}.
\end{equation*}
By repeating the same arguments as before, we obtain the following bound, which holds with probability $1-6\delta$:
\begin{align}
	U(f(X_P^*)) 
	&\leq U(f(X_n)) + P (U(f(X_n \cup \{\mathbf{x}^*\})) - U(f(X_n)))
	\label{eqn:proof:aeui_a}
	\\
	&\leq U(f(X_n)) + P \Bigl(
	\mathbb{E}_{p(\hat{f}|\mathcal{D}_n)}[U(\hat{f}(X_n \cup \{\mathbf{x}^*\}))] 
	+ \mathcal{A}_{\delta}(\mathbf{x}^*|\mathcal{D}_n) - U(f(X_n)) \Bigr) 
	\label{eqn:proof:aeui_b}
	\\
	&\leq U(f(X_n)) 
	+ P \Bigl(\frac{1}{H}\sum_{h=1}^H\hat{U}_J(\hat{f}_h(X_n \cup \{\mathbf{x}^*\})) 
	+ \mathcal{A}_\delta(\mathbf{x}^*| \mathcal{D}_n) - U(f(X_n)) + \xi_n \Bigr) 
	\label{eqn:proof:aeui_c}
	\\
	&\leq U(f(X_n)) + P \Bigl(\frac{1}{H}\sum_{h=1}^H\hat{U}_J(\hat{f}_h(X_{n+1})) + \mathcal{A}_\delta(\mathbf{x}_{n+1}| \mathcal{D}_n) - U(f(X_n)) + \xi_n \Bigr) 
	\label{eqn:proof:aeui_d}
	\\
	&\leq U(f(X_n)) + P (U(f(X_{n+1})) - U(f(X_n))) + 2P(\mathcal{A}_{\delta}(\mathbf{x}_{n+1}| \mathcal{D}_n) + \xi_n).
	\label{eqn:proof:aeui_e}
\end{align}
The first line \eqref{eqn:proof:aeui_a} follows from equation \eqref{eqn:submodular_bound_4}. The second line \eqref{eqn:proof:aeui_b} follows from applying the concentration bound \eqref{eqn:concentration_inequality} in \cref{ass:concentration}. The third line \eqref{eqn:proof:aeui_c} follows from applying Hoeffding's inequalities above. The fourth line \eqref{eqn:proof:aeui_d} is obtained by using the fact that the point $\mathbf{x}_{n+1} \in \mathbb{X}$ is picked according to the approximate AEUI acquisition function. The final line \eqref{eqn:proof:aeui_e} follows by applying the concentration inequalities again. As before, we can apply \cref{lemma:performance_bound} in order to obtain the final performance bound
\begin{equation*}
	U(f(X_N)) 
	\geq (1-e^{-N/P}) U(f(X_P^*)) - 2 \sum_{n=1}^N (\mathcal{A}_{\delta}(\mathbf{x}_n| \mathcal{D}_{n-1}) + \xi_{n-1} ) \Bigl(1-\frac{1}{P} \Bigr)^{N-n},
\end{equation*}
which holds with probability $1-6N\delta$. By using the inequality, $1-x \leq e^{-x}$ for all $x \in \mathbb{R}$, we obtain the final performance bound.

\begin{flushright}
	$\blacksquare$
\end{flushright}

\section{Additional figures} 
\label{app:additional_figures}
\begin{figure}[!h]
	\includegraphics[width=1\linewidth]{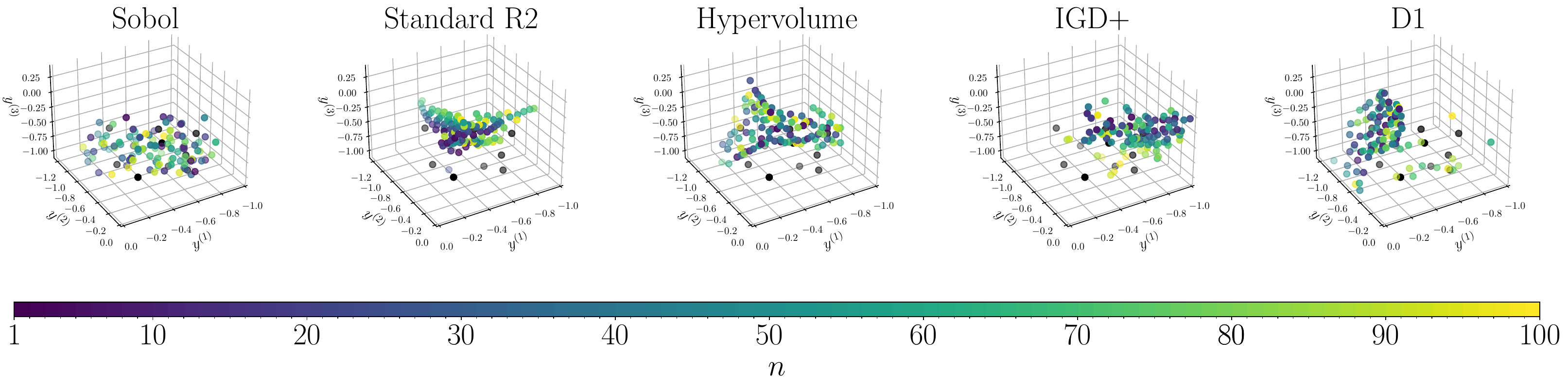}
	\includegraphics[width=1\linewidth]{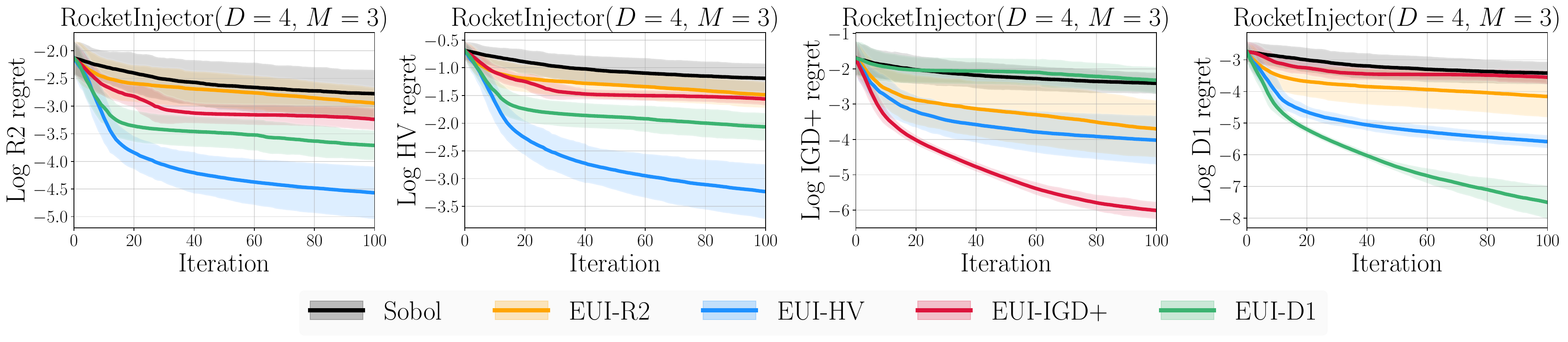}
	\centering
	\caption{An illustration of the EUI acquisition functions on the rocket injector problem. On the top row, we present a snapshot of the points obtained after one run of the Bayesian optimisation algorithm. On the bottom, we present the resulting performance plots.}
	\label{fig:rocket_performance}
\end{figure}

\begin{figure}[!h]
	\includegraphics[width=1\linewidth]{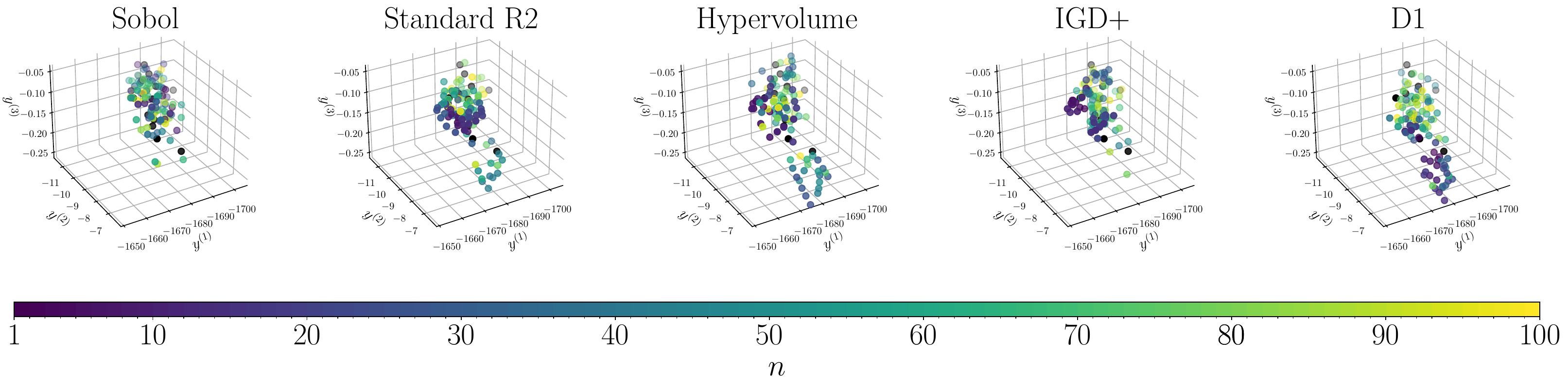}
	\includegraphics[width=1\linewidth]{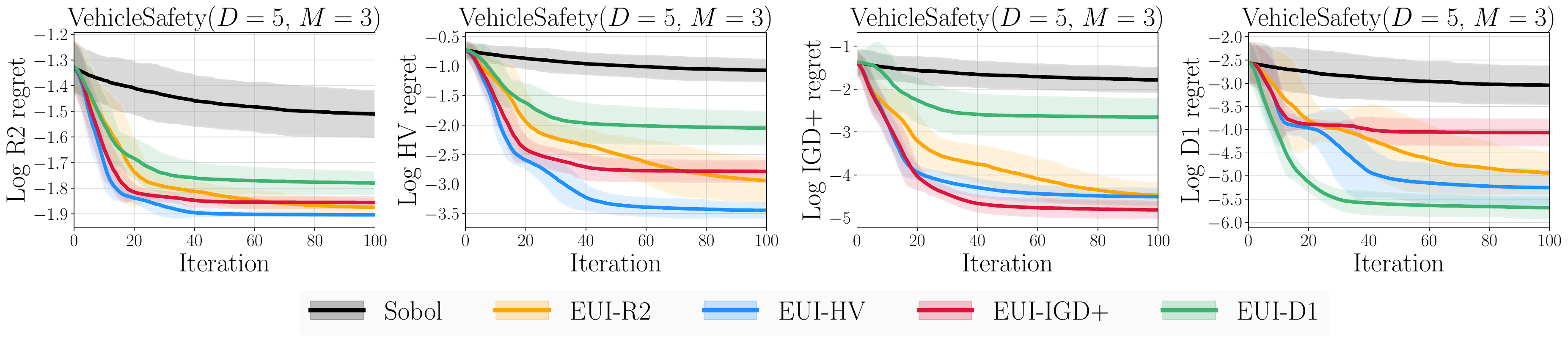}
	\centering
	\caption{An illustration of the EUI acquisition functions on the vehicle safety problem. On the top row, we present a snapshot of the points obtained after one run of the Bayesian optimisation algorithm. On the bottom, we present the resulting performance plots.}
	\label{fig:vehicle_performance}
\end{figure}

For some more intuition, we include two extra figures, \cref{fig:rocket_performance,fig:vehicle_performance}, which illustrates the impact of the choice of utility function when we use the Bayesian optimisation algorithm with the EUI acquisition function \eqref{eqn:eui}. In both of these examples, the number of objectives is $M=3$ and the utilities that we consider are the four special cases outlined in \cref{sec:special_cases}. For both of these problems, the standard R2 utility and hypervolume indicator are set to target the whole Pareto front, whilst the IGD+ and D1 utilities are used to target two complementary regions of the Pareto front. To elaborate, for the rocket injector problem (\cref{fig:rocket_performance}), we used the IGD+ utility to target the part of the Pareto front where the second objective is larger than the first objective in the normalised space, whilst the D1 utility is used to target the other part of the front. However, for the vehicle safety problem (\cref{fig:vehicle_performance}), we used the IGD+ utility to target the part of the Pareto front where the third objective is larger than the first objective in the normalised space, whilst the D1 utility is used to target the other part of the front. Overall, we find that the results here support the findings that we made earlier in \cref{sec:different_utilities}: namely, it is generally beneficial to optimise the same utility that is used to assess the performance, although there do exist some cases where using a different utility function leads to better performance. This latter result is likely attributed to problem-dependent features such as the geometry of the objective function and its Pareto front.

\newpage
\bibliographystyle{abbrv}
\bibliography{ms}
\end{document}